\documentclass{aims}

\usepackage{amsmath}
  \usepackage{paralist}
  \usepackage{graphics} 
  \usepackage{epsfig} 
\usepackage{lipsum}
\usepackage{amsfonts}
\usepackage{amsmath}
\usepackage{graphicx}
\usepackage{dsfont}
\usepackage{epstopdf}
\usepackage{algorithmic}
\usepackage{booktabs}
\usepackage{enumerate}
\usepackage{algorithm}
\usepackage{mathtools}
\usepackage[toc]{appendix}
\usepackage{multirow}
\usepackage{caption}
\usepackage{subcaption}
\usepackage{scalerel,stackengine}
\newcommand{\reallywidehat}[1]{%
\savestack{\tmpbox}{\stretchto{%
\scaleto{%
\scalerel*[\widthof{\ensuremath{#1}}]                         {\kern-.6pt\bigwedge\kern-.6pt}%
{\rule[-\textheight/2]{1ex}{\textheight}}
}{\textheight}%
}{0.5ex}}%
\ensurestackMath{\stackon[1pt]{#1}{\tmpbox}}%
}
\makeatletter
\newcommand*{\rom}[1]{\expandafter\@slowromancap\romannumeral #1@}
\makeatother
\ifpdf
  \DeclareGraphicsExtensions{.eps,.pdf,.png,.jpg}
\else
  \DeclareGraphicsExtensions{.eps}
\fi

\def\RR{\mathbb{R}}

\def\BU{\mathbf{U}}
\def\BV{\mathbf{V}}

\def\MPOD{{{\sc ho-pod}}}
\def\MODEIM{{{\sc ho-deim}}}
\def\MPDEIM{{{\sc ho-pod-deim}}}
\def\podeim{{{\sc pod-deim}}}

\usepackage{epstopdf}
 \usepackage[colorlinks=true]{hyperref}
\hypersetup{urlcolor=blue, citecolor=red}

\def\currentvolume{X}
 \def\currentissue{X}
  \def\currentyear{200X}
   \def\currentmonth{XX}
    \def\ppages{X--XX}

\newtheorem{example}{Example}

\newtheorem{proposition}{Proposition}

\theoremstyle{definition}

\newtheorem{remark}{Remark}

\title[Multilinear POD-DEIM for {semilinear} systems of ODEs] 
      {Multilinear POD-DEIM model reduction for 2D and 3D {semilinear} systems of differential equations}

\author[Gerhard Kirsten]{}


 \email{gerhard.kirsten2@unibo.it}

\thanks{The author is a member of Indam-GNCS, which support is gratefully acknowledged.}

\begin{document}
\maketitle

\centerline{\scshape Gerhard Kirsten$^*$}
\medskip
{\footnotesize
 \centerline{Dipartimento di Matematica, Universit{`a} di Bologna}
   \centerline{Piazza di Porta S. Donato, 5}
   \centerline{I-40127 Bologna, Italy}
} 

\bigskip


\begin{abstract}
We are interested in the numerical solution of coupled {\color{black}semilinear} partial differential equations (PDEs) in two and three dimensions. Under certain assumptions on the domain, we take advantage of the Kronecker structure arising in standard space discretizations of the differential operators and illustrate how the resulting system of ordinary differential equations (ODEs) can be treated directly in matrix or tensor form.  Moreover, in the framework of the proper orthogonal decomposition (POD) and the discrete empirical interpolation method (DEIM) we derive a two- and three-sided model order reduction strategy that is applied directly to the ODE system in matrix and tensor form respectively. We discuss how to integrate the reduced order model and, in particular, how to solve the tensor-valued linear system arising at each timestep of a semi-implicit time discretization scheme. We illustrate the efficiency of the proposed method through a comparison to existing techniques on classical benchmark problems such as the two- and three-dimensional Burgers equation.
\end{abstract}

\section{Introduction}

We are interested in the computationally efficient numerical solution of systems of {\color{black}semilinear} partial differential equations (PDEs) of the form
\begin{equation}\label{systempde}
\begin{cases}
\dot{u}_1 &= \mathcal{L}_{1}(u_1) + f_1(\nabla u_1, u_1,\ldots, u_{ g}, t)\\
&\vdots\\
\dot{u}_g &= \mathcal{L}_{ g}(u_{ g}) + f_{ g}(\nabla u_g , u_1,\ldots, u_{ g}, t),\\
\end{cases}
\end{equation}
where $u_i = u_i({\bf x}, t)$, with ${\bf x} \in \Omega \subset \RR^d$, $t \in [0,t_f]$ and suitable initial and boundary conditions, for all $i = 1,2,\ldots, g$. We restrict our attention to the two-dimensional and three-dimensional cases, that is $d = 2,3$. In this setting we assume that $\mathcal{ L}_i: \Omega\rightarrow \RR$ is linear in $u_i$, typically a diffusion operator, whereas $f_i: \Omega \times  [0,t_f] \rightarrow \RR$ is assumed to be nonlinear in $(\nabla u_i,u_1,\ldots u_{ g})$ and $t$. PDEs of the form \ref{systempde} describe mathematical models in several scientific fields, such as chemistry \cite{dewit99, vanag2004}, biology \cite{murray2001, hodgkin1952, george2013} and medicine \cite{sherratt2001}. For further applications we point the reader to \cite{Mainietal.01},%
\cite{Malchowetal.08},%
\cite{Quarteroni.17},%
\cite{Tveitoetal.10}, and references therein.

The simulation of \ref{systempde} and the accurate approximation to its solution pose several computational challenges. A semi-discretization in space of \ref{systempde} leads to a discrete approximation of the PDE inside a hypercube in $\RR^d$. The discrete model can then be integrated in time, most commonly by a time-discretization scheme, such as Implicit-Explicit schemes for instance \cite{Ascher1995,ruuth1995}. The method of lines (MOL) based on space discretizations, such as finite differences, finite elements, spectral methods, isogeometic analysis, rewrites the system \ref{systempde} as a system of ODEs of the form
\begin{equation}\label{systemode}
\begin{cases}
\dot{\bf u}_1(t) &= {\bf L}_1{\bf u}_1(t) + {\bf f}_1({\bf D}_1{\bf u}_1,{\bf u}_1,  \ldots, {\bf u}_{ g}, t),\\	
&\vdots\\
\dot{\bf u}_{ g}(t) &= {\bf L}_{ g}{\bf u}_{ g}(t) + {\bf f}_{ g}({\bf D}_g{\bf u}_g,{\bf u}_1, \ldots, {\bf u}_{ g}, t),\\	
\end{cases}
\end{equation}
where ${\bf u}_i \in \RR^{N}$ and each function ${\bf f}_i: \RR^N \times \RR^N \times \cdots \times \RR^N \times [0,t_f] \rightarrow \RR^N$ represents the function $f_i$ evaluated at the entries of the set of vectors $\{ {\bf u}_i\}_{i = 1}^{ g}$ and ${\bf D}_i{\bf u}_i$. In this setting, ${\bf L}_i \in \RR^{N \times N}$ accounts for the discretization of the linear operator $\mathcal{ L}_i$ on the selected basis and ${\bf D}_i \in \RR^{N \times N}$ for that of the gradient.  Furthermore, if finite differences are considered, then $N = \prod_{i = 1}^d n_i$, where $n_i$ represents the number of spatial nodes in the $x_i$ direction. It is clear that even for moderate $n_i$ the dimension of the considered matrices and vectors are very large when $d = 2,3$. Furthermore, in many cases fine grid discretizations of the system are required for an accurate simulation (see e.g., \cite{chat2011}). This poses a massive computational and memory challenge for systems of the form \ref{systemode}. 

 As a result, model order reduction techniques have been applied to dramatically reduce the dimension and the complexity of the resulting system of ODEs, see e.g., \cite{chat2011,wang2016, karasozen2020, sahyoun2014, karasozen2016, karasozen2021, kramer2011}. In particular, the common approach is to form a lexicographic ordering of the spatial nodes, unrolling the arrays (i.e., matrices or tensors when $d = 2$ and $d = 3$ respectively) of nodal values into long vectors in $\RR^N$, i.e., the unknown vectors ${\bf u}_i(t) \in \RR^N$. The dimension of the state space is then reduced to say $k \ll N$ through projection onto a low-dimensional subspace. Techniques for dimension reduction include  the proper orthogonal decomposition (POD)  \cite{benner2015review},\cite{Benner.05},\cite{hinze2005},\cite{kunisch1999},
reduced basis methods, see, e.g., \cite{patera2007reduced}, and rational interpolation strategies 
\cite{ABG.20} to mention a few. Furthermore, the complexity of the reduced model can be further decreased through hyper-reduction of the nonlinear term. This includes methods such as 
Missing Point Estimation (MPE) \cite{astrid2008}, the best points interpolation method (BPIM) \cite{nguyen2008}, and the discrete empirical interpolation method 
(DEIM) \cite{chaturantabut2010nonlinear}, which is based on the Empirical Interpolation Method 
(EIM), originally introduced in \cite{barrault2004}. Another model reduction approach is presented in \cite{kramer2019}. Here, the authors avoid the hyper-reduction step by rewriting the nonlinear term in polynomial form, through so-called lifting transformations; see e.g., \cite{gu2011}. The dimension of the resulting model is then fully reduced by POD.

 The success of the reduced order modelling is due to the fact that the solution trajectories are commonly attracted to low-dimensional manifolds \cite{chaturantabut2010nonlinear}. These low-dimensional models can then be rapidly simulated in a so-called online phase to form an approximation to the solution at the required timesteps. A shortcoming of the existing procedures, however, is the massive computational and storage demand in the offline phase. Even in the online phase, several vectors of length $N$ need to be stored in order to lift the low-dimensional functions back to the full dimension.  In this paper we aim to address precisely this shortcoming, with particular focus on  POD  for dimension reduction and DEIM for interpolation of the nonlinear function.

To this end, we illustrate that under certain hypotheses on the operators and the discretization basis in particular domains, \ref{systemode} can be equivalently expressed and integrated in array form, without forming a lexicographical ordering of the spatial nodes. In addition to a better structural interpretation of the discrete quantities, this can also lead to reduced memory requirements and computational costs. We discuss how to integrate the system in array form, and particularly also how to solve the tensor-structured linear system arising from the semi-implicit time integration when $d = 3$. 

Furthermore, inspired by \cite{kirsten2020}, we apply a $d-$sided POD-DEIM model order reduction, directly to the discrete system in array form, which is particularly advantageous, since no mapping is required from $\RR^{n_1 \times \cdots \times n_d}$ to $\RR^N$, and hence only vectors of length $n_1, \ldots, n_d$ need to be stored and processed. We mention that in the Arxiv report \cite[Section 8]{kirsten2020} the authors also present one example where a two-sided POD-DEIM order reduction strategy is applied to systems of matrix-valued ODEs. The presentation in this paper is, however, more general as it extends this idea to the higher dimensional multilinear array setting, including first order nonlinearities. Furthermore, a more in-depth experimental analysis of the two-dimensional case is presented than in \cite{kirsten2020}.

The paper is organized as follows. In section \ref{array} we illustrate how \ref{systemode} can be expressed in array form, whereas in section \ref{podreview} we review the standard {POD-DEIM} model reduction strategy. In section \ref{mpdeim} we extend POD-DEIM to the multilinear setting and illustrate how it can be applied to systems of array-valued ODEs in section \ref{sec:systemreduction}. The efficiency of the new procedure is illustrated by numerical experiments in section \ref{experiments} and our conclusions are formalized in section \ref{conclusion}.

{\it Notation.} 
{Scalar quantities are indicated by lower case letters and vectors are denoted by bold face lower case letters. Matrices are given by bold face upper case letters, whereas tensors are given by bold face, curly upper case letters and operators by standard curly upper case letters. ${\bf I}_n$ denotes the $n \times n$ identity matrix. For a matrix ${\bf M}$,
$\|{\bf M}\|$ denotes the matrix norm induced by the Euclidean vector norm, and 
$\|{\bf M}\|_F$ is the Frobenius norm. Furthermore, all reduced dimensional quantities are emphasized with a ` $\widehat{}$ '. 

For a third-order tensor $\pmb{ \mathcal{ T}} \in \RR^{n_1 \times n_2 \times n_3}$, the unfolding along the third mode is given by (see e.g., \cite{kolda2009})
$$
\pmb{ \mathcal{ T}}_{(3)} = \begin{pmatrix} {\bf T}_1, {\bf T}_2, \cdots, {\bf T}_{n_2} \end{pmatrix},$$
where $\pmb{ \mathcal{ T}}_{(3)}$ is a matrix in $\RR^{n_3 \times n_1n_2}$, and ${\bf T}_i \in \RR^{n_3 \times n_1}, i = 1,2,\ldots,n_2$
 is called a lateral slice. The multiplication of a tensor by a matrix, along a specific mode is done via the $m-$mode product, which, for a tensor $\pmb{ \mathcal{ T}} \in \RR^{n_1 \times n_2 \times n_3}$ and a matrix ${\bf M} \in \RR^{n \times n_m}$, we express as
$$
\pmb{\mathcal{ Q}} = \pmb{ \mathcal{ T}} \times_m {\bf M} \quad \iff \quad \pmb{\mathcal{ Q}}_{(m)} = {\bf M}\pmb{ \mathcal{ T}}_{(m)}.
$$
The Kronecker product of two matrices ${\bf M} \in \RR^{m_1 \times m_2}$ and ${\bf N} \in \RR^{n_1 \times n_2}$ is defined as
$$
{\bf M} \otimes {\bf N} = \begin{pmatrix}
M_{1,1}{\bf N} & \cdots &  M_{1,m_2}{\bf N}\\
\vdots & \ddots & \vdots\\
M_{m_1,1}{\bf N} & \cdots &  M_{m_2,m_2}{\bf N}
\end{pmatrix} \in \RR^{m_1n_1 \times m_2n_2},
$$
and the vec$(\cdot)$ operator maps the entries of a matrix, into a long vector, by stacking the columns of the matrix one after the other. The vectorization operator is applied to a third order tensor, via the first mode unfolding. Moreover, we will often make use of the property
\begin{equation}
\label{kronproperty}
({\bf M} \otimes {\bf N})\mbox{vec}({\bf X}) = \mbox{vec}({\bf N}{\bf X}{\bf M}^{\top}).
\end{equation}  
{\color{black}As a result, if $\pmb{\mathcal{ X}} \in \RR^{n_1 \times n_2 \times n_3}$, and ${\bf X} = \pmb{\mathcal{ X}}_{(3)}^{\top}$, then
\begin{equation}
\label{kronproperty2}
({\bf L} \otimes {\bf M} \otimes {\bf N})\mbox{vec}\left(\pmb{\mathcal{ X}}\right) = \mbox{vec}\left(({\bf M} \otimes {\bf N}){\bf X}{\bf L}^{\top}\right).
\end{equation}}
More properties used in the sequel are (see, e.g., \cite{golub13}): (i) $({\bf M} \otimes {\bf N})^{\top} = {\bf M}^{\top} \otimes {\bf N}^{\top}$; (ii) $ ({\bf M}_1 \otimes {\bf N}_1)({\bf M}_2 \otimes {\bf N}_2) = ({\bf M}_1{\bf M}_2 \otimes {\bf N}_1{\bf N}_2)$; (iii) $\| {\bf M} \otimes {\bf N} \|_2 = \|{\bf M}\|_2\|{\bf N}\|_2$; and (iv) $({\bf M} \otimes {\bf N})^{-1} = {\bf M}^{-1} \otimes {\bf N}^{-1}$,
where property (iv) holds if and only if both ${\bf M}$ and ${\bf N}$ are invertible.
}


All reported experiments were performed using MATLAB 9.9 (R2020b) 
(\cite{matlab2013}) on a MacBook Pro with 8-GB memory and a 2.3-GHz Intel core i5 processor.

\section{Matrix and tensor-based discretization of \ref{systempde}}
\label{array}
In this section we illustrate that under certain hypotheses, the discrete system \ref{systemode} can be expressed in terms of multilinear arrays; see e.g., \cite{Simoncini2017,Autilia2019matri, palitta2016}. To this end, suppose that $\mathcal{L}_1$,$\ldots$,$\mathcal{L}_g$ are second order differential operators with separable coefficients, such as the Laplace operator. Then, if $\mathcal{L}_i$ is discretized by means of a tensor basis, such as finite differences on parallelepipedal domains and certain spectral methods, the physical domain can be mapped to a reference hypercubic domain 
${\bf \Omega} = [a_1, b_1] \times \cdots \times [a_d, b_d]$. Hence, it holds that\footnote{We display the discretized Laplace operator, but more
general operators can also  be treated; see, e.g. \cite[Section 3]{Simoncini2017}.} (see e.g., \cite{palitta2016})
$$
{\bf L}_i = \sum_{m = 1}^{d} {\bf I}_{n_d} \otimes \cdots \otimes \overset{m}{\bf A}_{mi} \otimes \cdots \otimes {\bf I}_{n_1} \in \RR^{N \times N},
$$
and 
$$
{\bf D}_i = \sum_{m = 1}^{d} {\bf I}_{n_d} \otimes \cdots \otimes \overset{m}{\bf B}_{mi} \otimes \cdots \otimes {\bf I}_{n_1} \in \RR^{N \times N},
$$
where ${\bf A}_{mi} \in \RR^{n_m \times n_m}$ and ${\bf B}_{mi} \in \RR^{n_m \times n_m}$ contain the approximation of the second and first derivatives respectively in the $x_m$ direction, for $i = 1,2,\ldots,g$. The vectors ${\bf u}_i(t) \in \RR^N$ from \ref{systemode} then represent the vectorization of the elements of a tensor $\pmb{\mathcal{U}}_i(t) \in \RR^{n_1 \times \cdots \times n_d}$, such that ${\bf u}_i(t)  = {\tt vec}(\pmb{\mathcal{U}}_i(t))$, ${\bf L}_i{\bf u}_i = \mbox{vec}\left(\mathcal{A}_i(\pmb{\mathcal{U}}_i)\right)$ and ${\bf D}_i{\bf u}_i = \mbox{vec}\left(\mathcal{D}_i(\pmb{\mathcal{U}}_i)\right)$, where\footnote{For the case $d = 2$, \ref{A3} are Sylvester operators of the form ${\bf A}_{1i}{\bf U}_{i} + {\bf U}_{i}{\bf A}_{2i}^{\top}$ and ${\bf B}_{1i}{\bf U}_{i} + {\bf U}_{i}{\bf B}_{2i}^{\top}$ respectively \cite{Simoncini2017}.}
\begin{equation}\label{A3}
\mathcal{A}_i(\pmb{\mathcal{U}}_i)  := \sum_{m = 1}^d\, \pmb{\mathcal{U}}_i \times_m {\bf A}_{mi} \quad \mbox{and} \quad \mathcal{D}_i(\pmb{\mathcal{U}}_i)  := \sum_{m = 1}^d\, \pmb{\mathcal{U}}_i \times_m {\bf B}_{mi}.
\end{equation}
Moreover, if the function $\mathcal{F}_i: \RR^{n_1 \times \cdots \times n_d} \times \cdots \times  \RR^{n_1 \times \cdots \times n_d} \times [0,t_f] \rightarrow  \RR^{n_1 \times \cdots \times n_d}$ respresents the function $f_i$ evaluated at the entries of the arrays $\{\pmb{\mathcal{U}}_i\}_{i = 1}^g$ and $\mathcal{D}_i(\pmb{\mathcal{U}}_i)$, then it holds that ${\bf f}_i({\bf D}_i{\bf u}_i, {\bf u}_1, \ldots, {\bf u}_{ g},t) = {\tt vec}(\mathcal{F}_i(\mathcal{D}_i(\pmb{\mathcal{U}}_i),\pmb{\mathcal{U}}_1,\ldots, \pmb{\mathcal{U}}_g, t))$, and \ref{systemode} can be written in the form
\begin{equation}\label{arraybigsystem}
\begin{cases}
\dot{  \pmb{\mathcal{U}}  }_1 &= \mathcal{A}_1(\pmb{\mathcal{U}}_1) + \mathcal{F}_1\left(\mathcal{D}_1(\pmb{\mathcal{U}}_1),\pmb{\mathcal{U}}_1, \pmb{\mathcal{U}}_2, \cdots, \pmb{\mathcal{U}}_{ g}, t\right)\\
&\vdots\\
\dot{  \pmb{\mathcal{U}}  }_g &= \mathcal{A}_g(\pmb{\mathcal{U}}_g) + \mathcal{F}_g\left(\mathcal{D}_g(\pmb{\mathcal{U}}_g),\pmb{\mathcal{U}}_1, \pmb{\mathcal{U}}_2, \cdots, \pmb{\mathcal{U}}_{ g}, t\right),\\
\end{cases}
\end{equation}
with suitable initial conditions. The boundary conditions are contained in the matrices ${\bf A}_{mi}$ and ${\bf B}_{mi}$, $i = 1,\dots,g$, $m = 1,\ldots,d$; see e.g., \cite{Autilia2019matri, palitta2016}.

To simplify the presentation, we will consider the case where $n_1 = \cdots = n_d = n$  in the sequel, so that $N = n^d$. The extension to the more general case where $n_1 \ne \cdots \ne n_d$ is, however, possible.

\section{Review of POD-DEIM}
\label{podreview}
In this section we review the standard POD-DEIM method and its application to the 
dynamical system \ref{systemode}, with $g = 1$ and a nonlinear function independent of the gradient. We aim to reduce the dimension and complexity of a system of ODEs of the form 
\begin{equation}
\dot{\bf u}(t) = {\bf L}{\bf u}(t) + {\bf f}({\bf u},t), \qquad {\bf u}(0) = {\bf u}_0,
\end{equation}
with ${\bf u}(t) \in \RR^{N}$. Commonly POD is used to reduce the dimension of the state space by projecting onto a subspace of dimension $k \ll N$. In particular, given a set of $n_s$ trajectories (commonly known as snapshots) of the solution $\{{\bf u}(t_{j})\}_{j = 1}^{n_s}$, the POD basis is determined as the best rank $k$ approximation in the 2-norm of the space of snapshots. In particular, for 
$$
{\bf S} = [{\bf u}(t_1), \ldots, {\bf u}(t_{n_s})] \in \RR^{N \times n_s},
$$
the POD basis $\{{\bf v}_1,\ldots,{\bf v}_k\}$ is determined as the first $k$ dominant left singular vectors of ${\bf S}$. Therefore, if ${\bf V} = [{\bf v}_1,\ldots,{\bf v}_k] \in \RR^{N \times k}$, then we determine an approximation to ${\bf u}(t)$ as ${\bf u}(t) \approx {\bf V}\hat{\bf u}(t)$, where $\hat{\bf u}(t) \in \RR^{k}$ is determined as the solution of the reduced problem
\begin{equation}\label{redode}
\dot{\hat{\bf u}}(t) = \widehat{\bf L}\hat{\bf u}(t) + \widehat{\bf f}(\hat{\bf u},t), \qquad \hat{\bf u}(0) = {\bf V}^{\top}{\bf u}_0,
\end{equation}
with $\widehat{\bf L} = {\bf V}^{\top}{\bf L}{\bf V}$ and $\widehat{\bf f}(\hat{\bf u},t) = {\bf V}^{\top}{\bf f}({\bf V}\hat{\bf u},t)$. By definition, $\widehat{\bf f}(\hat{\bf u},t)$ first needs to be evaluated in full dimension, that is at the entries of ${\bf V}\hat{\bf u}(t) \in \RR^{N}$, before projection onto the low-dimensional subspace. Hence, the overall cost of evaluating \ref{redode} still depends on the full dimension $N$. One way to treat this bottleneck is through DEIM \cite{chaturantabut2010nonlinear}.

DEIM is used to interpolate a nonlinear function on the columns of an empirical basis. In particular if we consider a set of snapshots of the nonlinear function $\{{\bf f}({\bf u}, t_i)\}_{i = 1}^{n_s}$, then the DEIM basis $\{{\bf \phi}_1,\ldots,{\bf \phi}_p\}$ is determined as the first $p$ left singular vectors of the matrix
$$
{\bf N} = [{\bf f}({\bf u}, t_i), \ldots, {\bf f}({\bf u}, t_{n_s})] \in \RR^{N \times n_s}.
$$
Then, if we set ${\bf \Phi} = [{\bf \phi}_1,\ldots,{\bf \phi}_p] \in \RR^{N \times p} $ and consider ${\bf P} = [{\bf e}_{\rho_1},\ldots,{\bf e}_{\rho_p}] \in \RR^{N \times p}$ as a subset of columns of the $N\times N$ identity matrix, then the DEIM approximation of the nonlinear function is written as
\begin{equation}\label{deim}
{\bf f}({\bf u}, t) \approx {\bf \Phi}\left({\bf P}^{\top}{\bf \Phi}\right)^{-1}{\bf P}^{\top}{\bf f}({\bf u}, t).
\end{equation}
The interpolation indices $\{\rho_1, \ldots, \rho_p\}$ can be determined either in a greedy fashion \cite{chaturantabut2010nonlinear}, or through the pivoted $QR$ decomposition of ${\bf \Phi}^{\top}$ \cite{gugercin2018}. In the sequel we will make use of the latter algorithm and we will refer to it as {\tt q-deim}. The approximation \ref{deim} is particularly advantageous when the function ${\bf f}$ is evaluated elementwise at the entries of ${\bf u}$. In this case it holds that ${\bf P}^{\top}{\bf f}({\bf u}, t) = {\bf f}({\bf P}^{\top}{\bf u}, t)$, and hence ${\bf f}$ only needs to be evaluated at $p$ entries.

To complete the reduction of \ref{redode} to be independent of $N$, we therefore approximate $\widehat{\bf f}$ by \ref{deim}, so that
\begin{equation}
\widehat{\bf f}(\hat{\bf u}, t) \approx {\bf V}^{\top}{\bf \Phi}\left({\bf P}^{\top}{\bf \Phi}\right)^{-1}{\bf f}({\bf P}^{\top}{\bf V}_{k}\hat{\bf u}, t).
\end{equation}

In what follows we illustrate how the POD-DEIM method can be extended to the multilinear setting.

\section{POD-DEIM in the multilinear setting}
\label{mpdeim}
In this section we extend POD-DEIM to the matrix and tensor setting. We illustrate the procedure for systems of the form \ref{arraybigsystem} with $g = 1$ and a gradient-independent nonlinearity. The extension to the case of general $ g$ is presented later in the paper. In particular, we want to approximate the solution $\pmb{\mathcal{U}}(t) \in \RR^{n \times \cdots \times n}$, for $t \in [0,t_f]$, of the equation
\begin{equation}\label{tensorarray}
\dot{\pmb{\mathcal{U}}} = \mathcal{A}({\pmb{\mathcal{U}}}) + \mathcal{F}(\pmb{\mathcal{U}}, t), \qquad {\pmb{\mathcal{U}}}(0) = {\pmb{\mathcal{U}}}_0 ,
\end{equation}
by constructing $d$ low-dimensional basis matrices (one for each spatial mode $m = 1,\ldots,d$) ${\bf V}_m \in \RR^{n \times k_m}$, with $k_m \ll n$, to approximate $\pmb{\mathcal{U}}(t)$ in low dimension, for all $t \in [0,t_f]$.
To this end, given a set of snapshots $\{\pmb{\mathcal{U}}(t_{j})\}_{j = 1}^{n_s}$ with $\pmb{\mathcal{U}}(t_{j}) \in \RR^{n \times \cdots \times n}$, 
{\color{black} we define $\pmb{\mathcal{S}} \in \RR^{n \times \cdots \times n \times n_s}$ as a snapshot tensor\footnote{\color{black} We emphasize that this large dimensional tensor will never be explicitly formed or stored.} of order $d+1$ containing a collection of all the snapshots. In the well known case where the snapshots are vectors, this operation corresponds to collecting the vector snapshots into a \emph{snapshot matrix}. Instead, we are dealing with snapshots of higher dimension, which results in the definition of a \emph{snapshot tensor}. 

Then, given $\pmb{\mathcal{S}}$, }each matrix ${\bf V}_m$ is constructed in order to approximate the left range space of the matrix 
$$
\pmb{\mathcal{S}}_{(m)} = \begin{pmatrix} \pmb{\mathcal{U}}_{(m)}(t_1),\ldots,\pmb{\mathcal{U}}_{(m)}(t_{n_s})\end{pmatrix} \in \RR^{n \times n^dn_s}, \quad \mbox{for} \quad m = 1,\ldots,d,
$$
where $m$ represents the mode along which the tensor is unfolded. Forming or storing the matrix $\pmb{\mathcal{S}}_{(m)}$ is too computationally demanding,  even for moderate $n$ and $n_s$. Instead, the approximation spaces are updated one snapshot at a time. 

To this end, we determine the sequentially truncated higher order SVD\footnote{For the case $d = 2$, however, we just use the standard MATLAB {\tt SVD} function.} (STHOSVD) \cite{vannieuwenhoven2012} of each snapshot $\pmb{\mathcal{U}}(t_{j})$, so that
$
\pmb{\mathcal{U}}(t_{j}) \approx \pmb{\mathcal{C}}(t_{j}) \bigtimes_{m = 1}^d  \widetilde{\bf V}_{m}^{(j)},
$
where $ \widetilde{\bf V}_{m}^{(j)}$ contains the dominant left singular vectors of $\pmb{\mathcal{U}}_{(m)}(t_{j})$, truncated with respect to the tolerance $\tau$. {\color{black} Furthermore, $\pmb{\mathcal{C}}(t_{j})$ is the core tensor related to the STHOSVD of $\pmb{\mathcal{U}}(t_{j})$, and is generally defined as $\pmb{\mathcal{C}}(t_{j}) = \pmb{\mathcal{U}}(t_{j}) \bigtimes_{m = 1}^d  (\widetilde{\bf V}_{m}^{(j)})^{\top}$. Note, however, that in our procedure it is not necessary to explicitly compute the core tensors $\pmb{\mathcal{C}}(t_{j})$, since only the matrices $\widetilde{\bf V}_{m}^{(j)}$ are required.} More precisely, the approximation space in each mode is updated by orthogonalizing $\widetilde{\bf V}_{m}^{(j)}$ with respect to the previous basis vectors in that mode, and pruning any redundant information, with respect to $\tau$, to update ${\bf V}_m$. Through this procedure each snapshot can be discarded after it has been processed. We will refer to this as the higher order POD ({\sc ho-pod}) approximation. We call this a higher order approximation, since the standard POD projection is only a one-sided approximation, whereas this procedure is two-sided or three-sided for $d = 2$ and $d = 3$ respectively.

{\color{black} This type of approximation can also be interpreted as a $\mbox{Tucker}_{d}$ decomposition (see e.g., \cite[Section 4]{kolda2009}) of the snapshot tensor $\pmb{\mathcal{S}} \in \RR^{n \times \cdots \times n \times n_s}$ of the form 
$$
\pmb{\mathcal{S}} \approx \widehat{\pmb{\mathcal{S}}} \times_1 {\bf V}_{1} \times_2 \cdots \times_d {\bf V}_{d} \times_{d+1} {\bf I}_{n_s}, \quad  \widehat{\pmb{\mathcal{S}}} \in \RR^{k_1 \times \cdots \times k_d \times n_s},
$$
since the principal components of each of the first $d$ modes are analyzed. Nevertheless, instead of determining the core tensor $\widehat{\pmb{\mathcal{S}}}$, which contains a collection of low-dimensional approximations to all the given snapshots, we aim to use the basis matrices to approximate the solution of \ref{tensorarray} at time instances other than the ones considered for the snapshots. 

More precisely, we look for a approximation to the solution of \ref{tensorarray} of the form}
$
\pmb{\mathcal{U}}(t) \approx \widetilde{\pmb{\mathcal{U}}}(t) := \widehat{\pmb{\mathcal{U}}}(t) \bigtimes_{m=1}^d {\bf V}_m,
$
where  $\widehat{\pmb{\mathcal{U}}}(t) \in \RR^{k_1 \times \cdots \times k_d}$ $(k_m \ll n)$ satisfies the low-dimensional equation
\begin{equation}\label{tensorarraysmall}
\dot{\widehat{\pmb{\mathcal{U}}}} = \widehat{\mathcal{A}}(\widehat{\pmb{\mathcal{U}}}) + \widehat{\mathcal{F}}(\widehat{\pmb{\mathcal{U}}}, t), \qquad \widehat{\pmb{\mathcal{U}}}(0) = \widehat{\pmb{\mathcal{U}}}_0 ,
\end{equation}
%
%
where
\begin{equation}
\label{ABC}
\widehat{\mathcal{A}}(\widehat{\pmb{\mathcal{U}}})  := \sum_{m = 1}^d\, \widehat{\pmb{\mathcal{U}}} \times_m \widehat{\bf A}_{m}, \quad \widehat{\bf A}_m = {\bf V}_{m}^{\top}{\bf A}_m{\bf V}_{m},\quad \widehat{\pmb{\mathcal{U}}}(0) = \pmb{\mathcal{U}}_0 \bigtimes_{m = 1}^d {\bf V}_m^{\top}
\end{equation} and 
\begin{equation}\label{Flift}
\widehat{\mathcal{F}}(\widehat{\pmb{\mathcal{U}}}, t) = \mathcal{F}(\widetilde{\pmb{\mathcal{U}}}(t), t)  \bigtimes_{m = 1}^d {\bf V}_m^{\top}.
\end{equation}
{\color{black} For the time discretization of \ref{tensorarraysmall}, several alternatives can be considered, however it is well known that these type of semilinear equations are typically characterized by a stiff linear term and a nonstiff nonlinear term; see e.g., \cite{strikwerda2004}. That is, explicit methods will require unrealistically small timesteps to ensure stability in the linear term, whereas fully implicit schemes require the application of an expensive iterative nonlinear solver at each timestep. Instead, a good compromise is reached through semi-implicit (also called implicit-explicit) schemes, where the linear term is treated implicitly and the nonlinear term explicitly \cite[chapter IV.3]{hundsdorfer2013}. To this end,} we consider a second order implicit-explicit scheme, also known as IMEX 2--SBDF; see e.g., \cite{Autilia2019matri,Ascher1995}. Therefore, if $\widehat{\pmb{\mathcal{U}}}^{(j)}$ is an approximation of $\widehat{\pmb{\mathcal{U}}}(t_{j})$, then the linear system 
\begin{equation}\label{arraysmalltime}
(3\,\widehat{\mathcal{I}} - 2\Delta t\widehat{\mathcal{A}})(\widehat{\pmb{\mathcal{U}}}^{(j)}) = \widehat{\mathcal{G}}(\widehat{\pmb{\mathcal{U}}}^{(j-1)}, \widehat{\pmb{\mathcal{U}}}^{(j-2)})
\end{equation}
%
needs to be solved for each $t_{j}$, where 
\begin{equation*}\label{gdef}
\widehat{\mathcal{G}}(\widehat{\pmb{\mathcal{U}}}^{(j-1)}, \widehat{\pmb{\mathcal{U}}}^{(j-2)}) = 4\,\widehat{\pmb{\mathcal{U}}}^{({j-1})} -\,\, \widehat{\pmb{\mathcal{U}}}^{({j-2})} + 2\Delta t \left(2\widehat{\mathcal{F}}(\widehat{\pmb{\mathcal{U}}}^{({j-1})}, t_{j-1}) - \widehat{\mathcal{F}}(\widehat{\pmb{\mathcal{U}}}^{({j-2})}, t_{j-2})\right),
\end{equation*} 
and $\widehat{\mathcal{I}}: \RR^{k_1 \times \cdots \times k_d} \rightarrow \RR^{k_1 \times \cdots \times k_d}$ is the identity operator in the reduced dimension. To initiate the procedure, $\widehat{\pmb{\mathcal{U}}}^{({1})}$ can be determined by a semi-implicit Euler scheme from the known array $\widehat{\pmb{\mathcal{U}}}^{({0})}$. In what follows we discuss how \ref{arraysmalltime} is solved for $d = 2,3$.

\subsection{The solution of the linear system \ref{arraysmalltime}}
\label{systemsmall}
The solution of \ref{arraysmalltime} is not trivial, especially when $d = 3$, given that the matrices $\widehat{\bf A}_m$, for $m = 1,\ldots,d$ are necessarily dense due to the projection. When $d = 2$, it holds that $\widehat{\pmb{\mathcal{U}}}^{(j)} = \widehat{\bf  U}^{(j)} \in \RR^{k_1 \times k_2}$, and $\widehat{\mathcal{A}}$ is a Sylvester operator \cite{Simoncini2017}, so that \ref{arraysmalltime} is equivalent to the Sylvester equation (see e.g., \cite{Autilia2019matri})
$$
(3{\bf I}_{k_1} - 2\Delta t \widehat{\bf A}_1)\widehat{\bf  U}^{(j)} + \widehat{\bf  U}^{(j)}(-2\Delta t \widehat{\bf A}_2^{\top}) = \widehat{\mathcal{G}}(\widehat{\bf U}^{(j-1)}, \widehat{\bf U}^{(j-2)}) .
$$
Details on how to solve the Sylvester equation can be found in \cite{Simoncini2017}. In our experiments we make use of the built-in MATLAB function {\tt lyap}.

{\color{black} For the case $d = 3$, a direct method designed specifically for dense third order tensor linear systems has recently been introduced in \cite{simoncini2020} for tensors with a rank-one right hand side. Here we illustrate how this method can be applied to solve the linear system \ref{arraysmalltime}, accounting for a right hand side with rank greater than one.} To ease the readability, we drop the superscript $(j)$ for the description of the inner solver, when it is clear from the context.

{\color{black} By definition, when $d = 3$, the left hand side of \ref{arraysmalltime} can be vectorized as 
{\small $$
\left({\bf I}_{k_3} \otimes 3{\bf I}_{k_2} \otimes {\bf I}_{k_1} - {\bf I}_{k_3} \otimes {\bf I}_{k_2} \otimes 2\Delta t\widehat{\bf A}_1 -  {\bf I}_{k_3} \otimes 2\Delta t\widehat{\bf A}_2 \otimes {\bf I}_{k_1}- 2\Delta t\widehat{\bf A}_3 \otimes {\bf I}_{k_2} \otimes {\bf I}_{k_1}\right )\mbox{vec}\left(\widehat{\pmb{\mathcal{U}}}\right).
$$}Therefore, if we let $\widehat{\pmb{{{X}}}} = \widehat{\pmb{\mathcal{U}}}_{(3)}^{\top} \in \RR^{k_1k_2 \times k_3}$, then} by the use of property \ref{kronproperty2}, \ref{arraysmalltime} can be recast into the Sylvester equation
\begin{equation}\label{tensorsmallsylvester}
\left(3\, {\bf I}_{k_2} \otimes {\bf I}_{k_1} - {\bf I}_{k_2} \otimes 2\,\Delta t\widehat{\bf A}_1 -  2\,\Delta t\widehat{\bf A}_2 \otimes {\bf I}_{k_1}\right)\widehat{\pmb{{{X}}}}+ \widehat{\pmb{{{X}}}}\left(-2\,\Delta t\widehat{\bf A}_3^{\top}\right)  = \widehat{\bf G}. \end{equation}
Here $\widehat{\bf G} = \widehat{\mathcal{G}}(\widehat{\pmb{{{X}}}}^{({j-1})}, \widehat{\pmb{{{X}}}}^{({j-2})})  \in \RR^{k_1k_2 \times k_3}$.
Due to the large left dimension of this Sylvester equation, solving this directly is still not feasible. Instead, as is shown in \cite{simoncini2020}, it is possible to solve a sequence of much smaller Sylvester equations. To this end, let $\widehat{\bf A}_3^{\top} = {\bf Q}{\bf R}{\bf Q}^{\top}$ be the Schur decomposition of $\widehat{\bf A}_3^{\top}$. Then, if $\widehat{\pmb{{{Y}}}} = \widehat{\pmb{{{X}}}}{\bf Q}$, it holds that
\begin{equation}\label{tensorsmallsylvester2}
\left(3\,{\bf I}_{k_2} \otimes {\bf I}_{k_1} - {\bf I}_{k_2} \otimes 2\,\Delta t\widehat{\bf A}_1 -  2\,\Delta t\widehat{\bf A}_2 \otimes {\bf I}_{k_1}\right)\widehat{\pmb{{{Y}}}}+ \widehat{\pmb{{{Y}}}}\left(-2\,\Delta t{\bf R}\right) = \widehat{\bf G}{\bf Q},
\end{equation}
where ${\bf R} \in {\mathbb R}^{k_3 \times k_3}$ is block upper triangular. Therefore, by following the ideas of \cite{simoncini2020} and repeatedly using the property \ref{kronproperty}, the solution of \ref{arraysmalltime}, unfolded in the third mode, is given by
$$
\widehat{\pmb{\mathcal{U}}}_{(3)} = (\widehat{\bf Y}{\bf Q}^{\top})^{\top} = {\bf Q}\begin{pmatrix}\mbox{vec}({\bf Z}_1)^{\top}; \ldots; \mbox{vec}({\bf Z}_{k_3})^{\top} \end{pmatrix} \in \RR^{k_3 \times k_1k_2}.
$$
Here, ${\bf Z}_{h} \in \RR^{k_1 \times k_2}$ $(h = 1,2,\ldots,{k_3})$ solves the smaller Sylvester equation
$$
\left((3-2\,\Delta t{\bf R}_{h,h}){\bf I}_{k_1} - 2\,\Delta t\widehat{\bf A}_1\right){\bf Z}_{h} + {\bf Z}_{h}\left(-2\,\Delta t\widehat{\bf A}_2^{\top}\right) = {\bf H}_{h} +2\,\Delta t{\bf J}_{{h}-1},
$$
where $\widehat{\bf G}{\bf Q}  = {\bf H} =  \begin{pmatrix} \mbox{vec}({\bf H}_1), \ldots, \mbox{vec}({\bf H}_{k_3}) \end{pmatrix} \in \RR^{k_1k_2 \times {k_3}}$, and ${\bf J}_{h-1} \in \RR^{k_1 \times k_2}$ is the matricization of
\begin{equation}
\label{makeJ}
\mbox{vec}({\bf J}_{h-1}) = \begin{pmatrix} \mbox{vec}({\bf Z}_1), \ldots, \mbox{vec}({\bf Z}_{{h}-1}) \end{pmatrix}{\bf R}_{h,1:h-1}.
\end{equation}
In the special case where all coefficient matrices are symmetric and positive definite, the procedure can be even further accelerated; see \cite{simoncini2020} for further details. We refer to this inner solver as the {\sc t3-sylv} solver.
\subsection{Interpolation of the nonlinear function by {\sc ho-deim}}
\label{deimsec}
To determine the right-hand side $\widehat{\mathcal{G}}(\widehat{\pmb{\mathcal{U}}}^{(j-1)}, \widehat{\pmb{\mathcal{U}}}^{(j-2)})$ at each $t_{j}$, it is required to evaluate the nonlinear function in full dimension, as per the definition of $\widehat{\mathcal{F}}(\widehat{\pmb{\mathcal{U}}}, t)$. Instead, we interpolate the nonlinear function through a higher order version of DEIM. Consider the $d$ low-dimensional orthonormal matrices ${\bf \Phi}_{m} \in \RR^{n \times p_m}$,  with $p_m \ll n$,  determined as the output of {\sc ho-pod} of the set of nonlinear snapshots $\{\mathcal{F}(\pmb{\mathcal{U}}(t_{j}), t_{j})\}_{j = 1}^{n_s}$, for $m = 1,2,\ldots,d$. Furthermore, consider the $d$ selection matrices ${\bf P}_{m} \in \RR^{n \times p_m}$, given as the output of {\tt q-deim} with input ${\bf \Phi}_m^{\top}$, for $m = 1,2,\ldots,d$. The {\sc ho-deim} approximation of \ref{Flift} is then given by
\begin{equation}\label{tensordeim}
\widehat{\mathcal{F}}(\widehat{\pmb{\mathcal{U}}}, t) \approx \mathcal{F}(\widetilde{\pmb{\mathcal{U}}},t)\bigtimes_{m = 1}^d {\bf V}_{m}^{\top}{\bf \Phi}_{m}({\bf P}_{m}^{\top}{\bf \Phi}_{m})^{-1}{\bf P}_m^{\top}.
\end{equation}
If $\mathcal{F}$ is evaluated elementwise at the components of $\widetilde{\pmb{\mathcal{U}}}$, then it holds that
\begin{equation}\label{tensoreq}
\reallywidehat{\mathcal{F}(\widetilde{\pmb{\mathcal{U}}},t)} := \mathcal{F}(\widetilde{\pmb{\mathcal{U}}},t)\bigtimes_{m = 1}^d {\bf P}_{m}^{\top} =  \mathcal{F}(\widetilde{\pmb{\mathcal{U}}} \bigtimes_{m=1}^d {\bf P}_m^{\top},t).
\end{equation}
Notice that $\mathcal{F}$ is then evaluated at $p_1p_2\cdots p_d \ll n^d$ entries. {\color{black} Next we remark on a potential strategy for further reducing the online cost of the procedure. This is, however, not considered in our implementations.}

 {\color{black} \begin{remark}
 \label{deimremark}
{For certain nonlinear functions, the evaluation of $p_1p_2\cdots p_d$ entries online may not be feasible. One possibility that can be considered is to further approximate the {\sc ho-deim} reduced nonlinear function by a matrix-DEIM (MDEIM) type of interpolation (see e.g., \cite{bonomi2017,negri2015}). More precisely, if we consider the {\sc ho-deim} approximations \ref{tensordeim} and \ref{tensoreq} then the snapshot matrix
$
{\mathfrak N} = [{\mathfrak f}(t_1), \ldots, {\mathfrak f}(t_{n_s})] \in \RR^{p_{1}\cdots p_{d} \times n_s},
$
can be considered, where
$${\mathfrak f}(t_j) = \mbox{vec}\left( \reallywidehat{\mathcal{F}(\widetilde{\pmb{\mathcal{U}}},t_j)} \bigtimes_{m = 1}^d ({\bf P}_{m}^{\top}{\bf \Phi}_{m})^{-1} \right) \in \RR^{p_{1}\cdots p_{d}}.$$
If we define $\pmb{\mathcal{M}}_q \in \RR^{p_{1} \times \cdots \times p_{d}}$ as the tensorization of the $q$th column of ${\mathfrak M} \in \RR^{p_{1}\cdots p_{d} \times {\mathfrak p}}$ - the matrix of dominant left singular vectors of ${\mathfrak N}$ - then
\begin{equation*}
\begin{split}
\widehat{\mathcal{F}}(\widehat{\pmb{\mathcal{U}}}, t) &\approx \sum_{q = 1}^{\mathfrak p}c_q\pmb{\mathcal{M}}_q \bigtimes_{m = 1}^d {\bf V}_{m}^{\top}{\bf \Phi}_{m}, \qquad c_q = \left[\left({\mathfrak P}^{\top}{\mathfrak M}\right)^{-1}{\bf f}\left({\mathfrak P}^{\top}\mbox{vec}\left(\widetilde{\pmb{\mathcal{U}}} \bigtimes_{m=1}^d {\bf P}_m^{\top}\right),t\right) \right]_{q},\end{split}
\end{equation*}
where the columns of ${\mathfrak P} \in \RR^{p_1\cdots p_d \times {\mathfrak p}}$ are related to the {\tt q-deim} interpolation indices of the matrix ${\mathfrak M}$. Note that the nonlinear function is now evaluated at ${\mathfrak p} \ll p_1 \cdots p_d$ entries, even though no vectors of length $N$ need to be stored. This type of approximation has a couple of drawbacks, however. Firstly, to determine the snapshots ${\mathfrak f}(t_j)$, another (far cheaper) offline phase will be required. Moreover, potential structural properties such as symmetries, which are preserved by the approximation \ref{tensordeim}, may be destroyed by the vectorization; see also the discussion in the companion manuscript \cite[Section 4]{kirsten2020}.
}
\end{remark}}

We next provide an error bound for the {\sc ho-deim} approximation \ref{tensordeim}, 
where we recall that the matrices ${\mathbf \Phi}_{m}$, for $m = 1,2,\ldots,d$ all have orthonormal columns.
This bound is a direct extension to the array setting of \cite[Lemma 3.2]{chaturantabut2010nonlinear}. A similar result can be found in \cite{kirsten2020} for $d = 2$.
\begin{proposition}
\label{errprop}
{\color{black}Let ${\bf Q}_m = {\bf \Phi}_{m}({\bf P}_{m}^{\top}{\bf \Phi}_{m})^{-1}{\bf P}_m^{\top}$, and consider an arbitrary tensor $\pmb{\mathcal{F}} \in \RR^{n \times \cdots \times n}$, so that}
$$
\widetilde{\pmb{\mathcal{F}}} = \pmb{\mathcal{F}}\bigtimes_{m = 1}^d {\bf \Phi}_{m}({\bf P}_{m}^{\top}{\bf \Phi}_{m})^{-1}{\bf P}_m^{\top} = \pmb{\mathcal{F}}\bigtimes_{m = 1}^d {\bf Q}_m.
$$
Then,
\begin{equation}
\label{errbound}
\| \pmb{\mathcal{F}} - \widetilde{\pmb{\mathcal{F}}} \|_F \le {c_{1}}{c_{2}}\cdots c_d\, 
\|\pmb{\mathcal{F}} - \pmb{\mathcal{F}}\bigtimes_{m = 1}^d {\bf \Phi}_m{\bf \Phi}_m^{\top} \|_F
\end{equation}
where $ c_m = \left \| ({\bf P}_{m}^{\top}{\mathbf \Phi}_{m})^{-1} \right \|_2$, for $m = 1,\ldots,d$. 
\end{proposition}

\begin{proof}
Let ${\bf f} = \mbox{vec}(\pmb{\mathcal{F}}) \in \RR^N$. Then, by the properties of the Kronecker product
{\small \begin{equation*}
\begin{split}
\| \pmb{\mathcal{F}} - \widetilde{\pmb{\mathcal{F}}} \|_F &= \| \mbox{vec}(\pmb{\mathcal{F}}) - \mbox{vec}(\widetilde{\pmb{\mathcal{F}}}) \|_2 = \| {\bf f} - ( {\bf Q}_{d} \otimes \cdots \otimes {\bf Q}_{1}){\bf f} \|_2 \\ &= \left \|{\bf f} - ({\mathbf \Phi}_d \otimes \cdots \otimes {\mathbf \Phi}_1 )\left(({\bf P}_{d} \otimes \cdots \otimes {\mathbf P}_1 )^{\top}({\mathbf \Phi}_d \otimes \cdots \otimes {\mathbf \Phi}_1 )\right)^{-1}({\bf P}_{d} \otimes \cdots \otimes {\mathbf P}_1 )^{\top}{\bf f}\right\|_2
\end{split}
\end{equation*}}
Therefore, by \cite[Lemma 3.2]{chaturantabut2010nonlinear}, 
{\footnotesize
\begin{eqnarray*}
{\| \pmb{\mathcal{F}} - \widetilde{\pmb{\mathcal{F}}} \|_F} &\le&
{
     \left\|\left(({\bf P}_{d} \otimes \cdots \otimes {\mathbf P}_1 )^{\top}({\mathbf \Phi}_d \otimes \cdots \otimes {\mathbf \Phi}_1 )\right)^{-1}\right\|_2
\left\|{\bf f} - ({\mathbf \Phi}_d \otimes \cdots \otimes {\mathbf \Phi}_1 )({\mathbf \Phi}_d \otimes \cdots \otimes {\mathbf \Phi}_1 )^{\top}{\bf f}\right\|_2
}
\\&=& 
{
\left\|({\bf P}_{d}^{\top}{\mathbf \Phi}_{d})^{-1}\right\|_2
\cdots\left\|({\bf P}_{2}^{\top}{\mathbf \Phi}_{2})^{-1}\right\|_2 
\left\|({\bf P}_{1}^{\top}{\mathbf \Phi}_{1})^{-1}\right\|_2 
\left\|\pmb{\mathcal{F}} - \pmb{\mathcal{F}} \bigtimes_{m = 1}^d {\bf \Phi}_m{\bf \Phi}_m^{\top} \right\|_F.}\,\, 
\end{eqnarray*}}\end{proof}
The accuracy of the {\sc ho-deim} approximation therefore depends on the contraction coefficients $c_m$, which are minimized by the use of {\tt q-deim}; see, e.g., \cite{gugercin2018}. Furthermore it depends on the accuracy of the {\sc ho-pod} bases, given by the term $\left\|\pmb{\mathcal{F}} - \pmb{\mathcal{F}} \bigtimes_{m = 1}^d {\bf \Phi}_m{\bf \Phi}_m^{\top} \right\|_F$.

The full offline/online {\sc ho-pod-deim} reduction procedure for reducing tensor-valued ODEs is presented below in {{algorithm}  {\sc ho-pod-deim}} for the case $d = 3$.
 \begin{algorithm}
{{\bf Algorithm}  {\sc ho-pod-deim} for Tensor ODEs, $d = 3$}
\hrule
\vskip 0.1in
Given: {\color{black} Coefficient matrices of \ref{tensorarray} and function $\mathcal{F}: \RR^{n \times \cdots \times n} \times [0, t_f] \rightarrow \RR^{n \times \cdots \times n}$} 
\vskip 0.1in
\emph{Offline:}
\begin{enumerate}
\item {\color{black} For each $ j =1, 2, \ldots, n_{s}$
\begin{itemize}
\item[(i)] Iteratively update $\{\BV_{m}\}_{m = 1}^3$ and $\{{\bf \Phi}_{m}\}_{m=1}^3$, for the snapshots $\pmb{\mathcal{U}}(t_j)$ and 
$\mathcal{F}(\pmb{\mathcal{U}}(t_j), t_j)$ respectively as \ref{tensorarray} is integrated in time  and discard the snapshots ({\sc ho-pod});
\end{itemize}}
\item Compute $\widehat{\bf A}_m$, for $m = 1,2,3$ and $\widehat{\pmb{\mathcal{U}}}(0)$ from (\ref{ABC});
\item Determine $\{{\bf P}_{m}\}_{m = 1}^3$ using {\color{black}\tt q-deim} ({\sc ho-deim});
\vskip 0.05in
\item {Precompute $\{{\mathbf V}_{m}^{\top}{\mathbf \Phi}_{m}({\bf P}_{m}^{\top}
{\mathbf \Phi}_{m})^{-1}\}_{m=1}^3$ and $\{{\bf P}_{m}^{\top}{\mathbf V}_{m}\}_{m=1}^3$} ;
\item Compute the Schur decomposition $\widehat{\bf A}_3^{\top} = {\bf Q}{\bf R}{\bf Q}^{\top}$;\\
\end{enumerate}
\vskip 0.1in
\emph{Online:}
\begin{enumerate}
\item Determine $\widehat{\pmb{\mathcal{U}}}^{(1)}$ from $\widehat{\pmb{\mathcal{U}}}^{(0)}$;
\item For each $ j =2, 3, \ldots, n_{\mathfrak t}$
\begin{itemize}
{\item[(i)] Approximate $\widehat{\mathcal{F}}(\widehat{\pmb{\mathcal{U}}}^{(j-1)}, t_{j - 1})$ and $\widehat{\mathcal{F}}(\widehat{\pmb{\mathcal{U}}}^{(j-2)}, t_{j - 2})$ as in \ref{tensordeim} and \ref{tensoreq}
using the matrices computed above, and evaluate ${\mathcal{G}}(\widehat{\pmb{\mathcal{U}}}^{(j-1)}, \widehat{\pmb{\mathcal{U}}}^{(j-2)})$;}
\item[(ii)] For each $ h =1, 2, \ldots, k_3$:
\begin{enumerate}
\item Evaluate ${\bf J}_{h - 1}$ using \ref{makeJ} and compute ${\bf H} = \widehat{\bf G}{\bf Q}$;
\item Reshape column $h$ of  ${\bf H}$ into a $k_1 \times k_2$ matrix to form ${\bf H}_{h}$ ;
\item Solve the Sylvester matrix equation by a direct solver:     
$$
\left((3 - 2\Delta t{\bf R}_{h,h}){\bf I}_k - 2\,\Delta t\widehat{\bf A}\right){\bf Z}_{h} + {\bf Z}_{h}\left(-2\,\Delta t\widehat{\bf B}^{\top}\right) = {\bf H}_{h} +2\,\Delta t{\bf J}_{{h}-1},
$$
\item Update ${\bf Z} \leftarrow [{\bf Z}, \mbox{vec}({\bf Z}_{h})]$;
\end{enumerate}
\item[(iii)] Evaluate $\widehat{\pmb{\mathcal{U}}}_{(3)}^{(j)} = {\bf Q}{\bf Z}^{\top}$ and reshape it into a $k_1 \times k_2 \times k_3$ tensor;
\end{itemize}
\item Return $\BV_{1}, \BV_{2}, \BV_{3}$ and $\{\widehat{\pmb{\mathcal{U}}}^{(j)}\}_{j = 1}^{n_{\mathfrak t}}$, so that $\widehat{\pmb{\mathcal{U}}}^{(j)} \times_1 \BV_{1} \times_2 \BV_{2} \times_3 \BV_{3} \approx \pmb{\mathcal{U}}(t_{j})$;
\end{enumerate}
\end{algorithm}
In what follows, we illustrate how the discussed higher-order POD-DEIM order reduction strategy can be applied to systems of ODEs of the form \ref{arraybigsystem}.

\section{Order reduction of systems of array-valued ODEs}
\label{sec:systemreduction}
Here we  illustrate how the {\sc ho-pod-deim} order reduction scheme presented in the previous section can be applied to systems of array-valued ODEs of the form \ref{arraybigsystem}. Indeed, consider $d\cdot g$ tall basis matrices ${\mathbf V}_{m,i} \in \RR^{n \times k_{mi}}$ with orthonormal columns, for $i = 1,2,\ldots,  g$
and $m = 1,2,\ldots,d$, where $k_{mi} \ll n$. {\color{black} That is, we consider $d$ basis matrices for each of the $g$ equations in \ref{arraybigsystem}.} Approximations to  each ${\pmb{\mathcal{U}}}_i(t)$, for $t \in  [0,t_f]$, can then be written as
$$
{\pmb{\mathcal{U}}}_i(t) \approx \widetilde{\pmb{\mathcal{U}}}_i(t) = \widehat{\pmb{\mathcal{U}}}_i(t) \bigtimes_{m =1}^d {\bf V}_{m,i}, \qquad i = 1,2,\ldots, g.
$$
The functions $\widehat{\pmb{\mathcal{U}}}_i(t) \in \RR^{k_{1i} \times k_{2i} \times \cdots \times k_{di}}$ are determined as an approximation to the solution of the reduced, coupled problem
\begin{equation}\label{arraysmallsystem}
\begin{cases}
\dot{ \widehat{ \pmb{ \mathcal{U} } } }_1 &= \widehat{\mathcal{A}}_1(\widehat{\pmb{\mathcal{U}}}_1) + \widehat{\mathcal{F}}_1\left(\widehat{\mathcal{D}}_1(\widehat{\pmb{\mathcal{U}}}_1),\widehat{\pmb{\mathcal{U}}}_1, \widehat{\pmb{\mathcal{U}}}_2, \cdots, \widehat{\pmb{\mathcal{U}}}_{ g}, t\right)\\
&\vdots\\
\dot{\widehat{\pmb{\mathcal{U}}}}_g &= \widehat{\mathcal{A}}_g(\widehat{\pmb{\mathcal{U}}}_g) + \widehat{\mathcal{F}}_g\left(\widehat{\mathcal{D}}_g(\widehat{\pmb{\mathcal{U}}}_g),\widehat{\pmb{\mathcal{U}}}_1, \widehat{\pmb{\mathcal{U}}}_2, \cdots, \widehat{\pmb{\mathcal{U}}}_{ g}, t\right),\\
\end{cases}
\end{equation}
where
\begin{equation}
\label{sysABC}
\widehat{\mathcal{A}}_i(\widehat{\pmb{\mathcal{U}}}_i)  := \sum_{m = 1}^d\, \widehat{\pmb{\mathcal{U}}}_i \times_m \widehat{\bf A}_{mi}, \quad \widehat{\bf A}_{mi} = {\bf V}_{m,i}^{\top}{\bf A}_{mi}{\bf V}_{m,i},\quad \widehat{\pmb{\mathcal{U}}}(0) = {\pmb{\mathcal{U}}}_{i0} \bigtimes_{m = 1}^d {\bf V}_{m,i}^{\top}
\end{equation} and 
\begin{equation}\label{sysFlift}
\widehat{\mathcal{F}}_i\left(\widehat{\mathcal{D}}_i(\widehat{\pmb{\mathcal{U}}}_i),\widehat{\pmb{\mathcal{U}}}_1, \cdots, \widehat{\pmb{\mathcal{U}}}_{ g}, t\right) = \mathcal{F}_i\left(\mathcal{D}_i(\widetilde{\pmb{\mathcal{U}}}_i),\widetilde{\pmb{\mathcal{U}}}_1, \cdots, \widetilde{\pmb{\mathcal{U}}}_{ g}, t\right)  \bigtimes_{m = 1}^d {\bf V}_{m,i}^{\top}.
\end{equation}
We use the {\sc ho-pod} procedure from the previous section to determine the basis matrices. In particular, given the set of snapshot solutions $\{ {\pmb{\mathcal{U}}}_i (t_{j})\}_{j = 1}^{n_s}$, the basis matrices ${\mathbf V}_{m,i} \in \RR^{n \times k_{mi}}$, $m = 1,2,\ldots,d$, are determined following the \MPOD\ procedure from section \ref{mpdeim}, for each $i = 1,2, \ldots,  g$.

The reduced order model \ref{arraysmallsystem}, can also be integrated by means of the IMEX 2~-~SBDF scheme for systems. Indeed, the $\widehat{\pmb{\mathcal{U}}}_{1}^{(j)}, \widehat{\pmb{\mathcal{U}}}_{2}^{(j)}, \ldots, \widehat{\pmb{\mathcal{U}}}_{ g}^{(j)}$ 
approximations to $\widehat{\pmb{\mathcal{U}}}_{1}(t_{j}), \widehat{\pmb{\mathcal{U}}}_{2}(t_{j}), \ldots, \widehat{\pmb{\mathcal{U}}}_{ g}(t_{j})$ are determined by solving the linear systems
\begin{equation}\label{arraysystemtime}
\begin{cases}
(3\,\widehat{\mathcal{I}} - 2\,\Delta t\widehat{\mathcal{A}}_1)(\widehat{\pmb{\mathcal{U}}}_1^{(j)}) &= \widehat{\mathcal{G}}_1\left(\{\widehat{\pmb{\mathcal{U}}}_i^{(j-1)}\}_{i = 1}^g, \{\widehat{\pmb{\mathcal{U}}}_i^{(j-2)}\}_{i = 1}^{g}\right)\\
&\vdots\\
(3\,\widehat{\mathcal{I}} - 2\,\Delta t\widehat{\mathcal{A}}_g)(\widehat{\pmb{\mathcal{U}}}_g^{(j)}) &= \widehat{\mathcal{G}}_g\left(\{\widehat{\pmb{\mathcal{U}}}_i^{(j-1)}\}_{i = 1}^g, \{\widehat{\pmb{\mathcal{U}}}_i^{(j-2)}\}_{i = 1}^{g}\right),\\
\end{cases}
\end{equation}
at each $t_j$, where
{\small\begin{equation*}\label{gdef2}
\begin{split}
\widehat{\mathcal{G}}_i\left(\{\widehat{\pmb{\mathcal{U}}}_i^{(j-1)}\}_{i = 1}^g, \{\widehat{\pmb{\mathcal{U}}}_i^{(j-2)}\}_{i = 1}^{g}\right) =  4\,\widehat{\pmb{\mathcal{U}}}_i^{({j-1})}  
&+ 4\Delta t  \widehat{\mathcal{F}}_i\left(\widehat{\mathcal{D}}_i(\widehat{\pmb{\mathcal{U}}}_i^{(j-1)}),\widehat{\pmb{\mathcal{U}}}_1^{(j-1)}, \cdots, \widehat{\pmb{\mathcal{U}}}_{ g}^{(j-1)}, t_{j-1}\right)\\ - \widehat{\pmb{\mathcal{U}}}_i^{({j-2})} &- 2\Delta t \widehat{\mathcal{F}}_i\left(\widehat{\mathcal{D}}_i(\widehat{\pmb{\mathcal{U}}}_i^{(j-2)}),\widehat{\pmb{\mathcal{U}}}_1^{(j-2)}, \cdots, \widehat{\pmb{\mathcal{U}}}_{ g}^{(j-2)}, t_{j-2}\right).
\end{split}
\end{equation*}}Each of the $g$ linear systems in \ref{arraysystemtime} can be solved via the procedures set forth in section \ref{systemsmall} for $d = 2,3$.

Adjacent to the setting discussed in section \ref{deimsec}, the nonlinear functions need to be interpolated to avoid evaluating the functions in full dimension at each timestep. To this end,
we approximate $\widehat{\mathcal{F}}_i\left(\widehat{\mathcal{D}}_i(\widehat{\pmb{\mathcal{U}}}_i),\widehat{\pmb{\mathcal{U}}}_1, \cdots, \widehat{\pmb{\mathcal{U}}}_{ g}, t\right) $ in the space spanned by the columns of the matrices ${\mathbf \Phi}_{m,i} \in \RR^{n \times p_{mi}}$, $m = 1,2,\ldots,d$, for each $i = 1,2,\ldots g$, where $p_{mi} \ll n$. Given the selection matrices ${\bf P}_{m, i} \in \RR^{n \times p_{mi}}$, $m = 1,2,\ldots,d$, we obtain
\begin{equation}
\begin{split}
\widehat{\mathcal{F}}_i(\widehat{\mathcal{D}}_i(\widehat{\pmb{\mathcal{U}}}_i), \widehat{\pmb{\mathcal{U}}}_1, \ldots, \widehat{\pmb{\mathcal{U}}}_{ g}, t) &\approx \mathcal{F}_i(\mathcal{D}_i(\widetilde{\pmb{\mathcal{U}}}_i),\widetilde{\pmb{\mathcal{U}}}_1,\ldots, \widetilde{\pmb{\mathcal{U}}}_{ g}, t)  \bigtimes_{m = 1}^d {\bf V}_{m,i}^{\top}{\bf Q}_{m,i},\\
\end{split}
\end{equation}
with the oblique projectors
$$
{\bf Q}_{m,i} = {\bf \Phi}_{m,i}({\bf P}_{m,i}^{\top}{\bf \Phi}_{m,i})^{-1}{\bf P}_{m,i}^{\top}, \quad \mbox{for} \quad m = 1,2,\ldots,d \quad \mbox{and} \quad i = 1,2,\ldots,  g.
$$
The basis matrices ${\mathbf \Phi}_{m,i}$, $m = 1,2,\ldots,d$ are determined via the {\sc ho-pod} procedure described in section \ref{mpdeim}, given the set of nonlinear snapshots 
\begin{equation}\label{nonsnaps}
\left \{\mathcal{F}_i\left(\mathcal{D}_i({\pmb{\mathcal{U}}}_i(t_j)),{\pmb{\mathcal{U}}}_1(t_j), {\pmb{\mathcal{U}}}_2(t_j), \cdots, {\pmb{\mathcal{U}}}_{ g}(t_j), t_j\right)\right\}_{j = 1}^{n_s},
\end{equation}
 whereas the selection matrices ${\mathbf P}_{m,i}$ are determined via {\tt q-deim} with inputs ${\mathbf \Phi}_{m,i}^{\top}$ respectively, for each $i = 1,2,\ldots,g$. In this paper we assume that there is a componentwise relationship between the arrays ${\pmb{\mathcal{U}}}_1(t_j), {\pmb{\mathcal{U}}}_2(t_j), \cdots, {\pmb{\mathcal{U}}}_{ g}(t_j)$ and the approximation to the gradient $\mathcal{D}_i({\pmb{\mathcal{U}}}_i)$ in the nonlinear function $\mathcal{F}_i$. Therefore, since the matrix ${\mathbf P}_{m,i}$ is merely responsible for selecting rows in the respective modes, it holds that
\begin{equation*}\label{applyp}
\begin{split}
\reallywidehat{\mathcal{F}_i\left(\mathcal{D}_i(\widetilde{\pmb{\mathcal{U}}}_i),\widetilde{\pmb{\mathcal{U}}}_1,\ldots, \widetilde{\pmb{\mathcal{U}}}_{ g}, t\right)} &:=  \mathcal{F}_i\left(\mathcal{D}_i(\widetilde{\pmb{\mathcal{U}}}_i),\widetilde{\pmb{\mathcal{U}}}_1,\ldots, \widetilde{\pmb{\mathcal{U}}}_{ g}, t\right)  \bigtimes_{m = 1}^d {\bf P}_{m,i}^{\top}\\
&= \mathcal{F}_i\left(\mathcal{D}_i(\widetilde{\pmb{\mathcal{U}}}_i) \bigtimes_{m = 1}^d {\bf P}_{m,i}^{\top},\,\,\widetilde{\pmb{\mathcal{U}}}_1 \bigtimes_{m = 1}^d {\bf P}_{m,i}^{\top},\ldots,\widetilde{\pmb{\mathcal{U}}}_{ g} \bigtimes_{m = 1}^d {\bf P}_{m,i}^{\top}, t\right).
\end{split}
\end{equation*}

\section{Numerical experiments}
\label{experiments}
In this section we illustrate the efficiency of the discussed methods via benchmark problems from biology and engineering. For all problems, the accuracy of the reduced order model is tested through the average error measure
\begin{equation}\label{errormeasure2}
\bar{\mathcal{E}}(\pmb{\mathcal{ U}})= \frac{1}{n_\mathfrak{t}}
\sum_{j = 1}^{n_\mathfrak{t}}\frac{\| \pmb{\mathcal{ U}}^{(j)}- \widetilde{\pmb{\mathcal{ U}}}^{(j)} \|_F}{\|\pmb{\mathcal{ U}}^{(j)}\|_F},
\end{equation}
and the truncation of the singular values is done by monitoring the quality of the approximation in the Frobenius norm. That is, if $\sigma_1 \ge \sigma_2 \ge \cdots \ge \sigma_{\kappa}$ are the singular values of the matrix that needs to be truncated, then the new dimension $\nu \le \kappa$ is determined as
\begin{equation} \label{kselect}
\frac{\sqrt{\sum_{i = \nu + 1}^{\kappa} {\sigma}_i^2}}{\sqrt{\sum_{i = 1}^{\kappa} {\sigma}_i^2}} < \tau.
\end{equation}
We first illustrate the efficiency with two examples where $d =2$, after which we investigate the three-dimensional coupled Burgers equation.
\begin{example}\label{ex:system}
{\rm {\it The 2D FitzHugh-Nagumo model (FN).} Consider the following classical problem, given in adimensional form,
\begin{equation}
\label{FHN}
\dot{u}_1 =  \delta_1\Delta u_1 + \Gamma (-u_1^3 + u_1 - u_2), \qquad 
\dot{u}_2 =  \delta_2\Delta u_2+ \Gamma (\beta u_1 - \beta h u_2) ,
\end{equation}
where the functions $u_1(x,y,t)$ and $u_2(x,y, t)$ model the densities of two species for $t \in [0,1]$, and $[x, y] \in  [-1, 1]^2$.
We refer the reader to, e.g., \cite{gambino2019} for a description of the role of the nonnegative 
coefficients $h, \beta,$ and $\Gamma$. 
For this example we set $h = 0.5, \beta = 2.1, \Gamma = 9.65$, $\delta_1 = 0.01$ and $\delta_2 = 0.1$. 
Furthermore, homogeneous Neumann boundary conditions are imposed and the initial state is given by
\begin{equation*}
\begin{split}
u_1(x, y, 0) &= (1-x^2)(1-y^2)\sin(2\pi x)\cos(2\pi (y+0.3)) \\
u_2(x, y, 0) &= (1-x^2)(1-y^2)e^{-\sin(2\pi(x-0.3)y)}.
\end{split}
\end{equation*}
This example investigates the efficiency of the reduced order model in terms of accuracy and online CPU time. To this end, the system \ref{FHN} is discretized with $n = 1200$ spatial nodes in each direction yielding
the form \ref{arraybigsystem}. Note that this is equivalent to the system \ref{systemode} with dimension $N = 1\,440\,000$.

 {\color{black} In particular, if we let ${\bf T} = \mbox{tridiag}(1,\underline{-2},1) + {\bf N}$, ${\bf T} \in \RR^{n \times n}$, where
$$
{\bf N} = \frac{2}{3}\begin{pmatrix} 2 & -1/2 & \cdots & 0 & 0 &\\
0 & 0 & \cdots & \cdots & 0 &\\
\vdots & & & \vdots &\\
0& 0& \cdots & -1/2 & 2 &
\end{pmatrix} \in \RR^{n \times n}
$$
contains the Neumann boundary conditions (see e.g., \cite{Autilia2019matri}), then the coefficient matrices of \ref{arraybigsystem} are defined as
$$
{\bf A}_{11} = \Gamma\, {\bf I}_n + \frac{\delta_1}{\ell_x^2}{\bf T}, \quad {\bf A}_{21} = -\Gamma \beta h{\bf I}_n + \frac{\delta_2}{\ell_x^2}{\bf T}, \quad {\bf A}_{12} =  \frac{\delta_1}{\ell_y^2}{\bf T},\quad \mbox{and} \quad {\bf A}_{22} =  \frac{\delta_2}{\ell_y^2}{\bf T}
$$
where $\ell_x = \ell_y = 2/(n-1)$. Notice that the discretized linear terms $\Gamma{\bf U}_1$ and $-\Gamma\beta h{\bf U}_2$ have been incorporated into the coefficient matrices ${\bf A}_{11}$ and ${\bf A}_{21}$ respecctively. Furthermore the matrix ${\mathcal{F}}_1(\BU_1,{\bf U}_2,t)$ stems from evaluating 
the function $f_1(u_1,u_2) = -\Gamma u_1^3$ elementwise, whereas ${\mathcal{F}}_2(\BU_1,{\bf U}_2,t) = \Gamma \beta {\bf U}_1$ is linear and requires no {\sc ho-deim} interpolation. Finally, the remaining linear term in the first equation $-\Gamma \, {\bf U}_2$ can be projected explicitly onto the {\sc ho-pod} subspace of the first equation.}

 In our experiments we found that $n_s = 20$ equispaced snapshots $\BU_1(t)$ and ${\bf U}_2(t)$ in the timespan $[0,1]$ are sufficient for constructing the basis vectors. Furthermore, we consider four different truncation tolerances ($\tau = 10^{-2}, 10^{-4}, 10^{-6}, 10^{-8}$) for this experiment. Table \ref{comp1} reports all the basis dimensions obtained for each $\tau$, by means of \ref{kselect}. 
\begin{table}[htb!]
\caption{\color{black}Example~\ref{ex:system}. Dim. of \MPOD\ and \MODEIM\ bases obtained for different $\tau$. The full order model has dimension $n = 1200$.  \label{comp1}}
\centering
\begin{tabular}{|l|l|c|c|c|c|}
\hline             
   &      &   left dim.       &  right dim.   & left dim.     & right dim.  \\
   $\tau$ & ${\bf U}_i$     &   \MPOD\          &  \MPOD\       & \MODEIM\       & \MODEIM\      \\ \hline
   \multirow{2}{*}{{$10^{-2}$}} & $\BU_1$  & 7         & 7  & 11       & 11         \\ 
 & $\BU_2$  & 9          & 10  & --        & --        \\ \hline
\multirow{2}{*}{{$10^{-4}$}} & $\BU_1$  & 18         & 20   & 23       & 23            \\ 
 & $\BU_2$  & 19         & 20  & --        & --          \\ \hline
 \multirow{2}{*}{{$10^{-6}$}} & $\BU_1$  & 31         & 33   & 32        & 34           \\ 
 & $\BU_2$  & 29       & 31  & --        & --         \\ \hline
  \multirow{2}{*}{{$10^{-8}$}} & $\BU_1$  & 43          & 46   & 44        & 47        \\ 
 & $\BU_2$  & 37       & 40  & --        & --            \\ \hline
\end{tabular}
\end{table}

In Figure \ref{fhnfig} (left) we plot the average error \ref{errormeasure2} for both $\widetilde{\BU}_1$ and $\widetilde{\BU}_2$ integrated from $0$ to $t_f$ at $n_\mathfrak{t} = 300$ timesteps, for the different values of $\tau$ presented in Table~\ref{comp1}. For the error computation, both the full order model and the reduced order model are integrated with the IMEX 2-SBDF scheme. On the right of Figure \ref{fhnfig} we plot the CPU time for integrating the full order model and the reduced order model at $n_\mathfrak{t} = 300$ timesteps for decreasing $\tau$. 
\begin{figure}[htb]
\minipage{0.49\textwidth}
  \includegraphics[width=\linewidth]{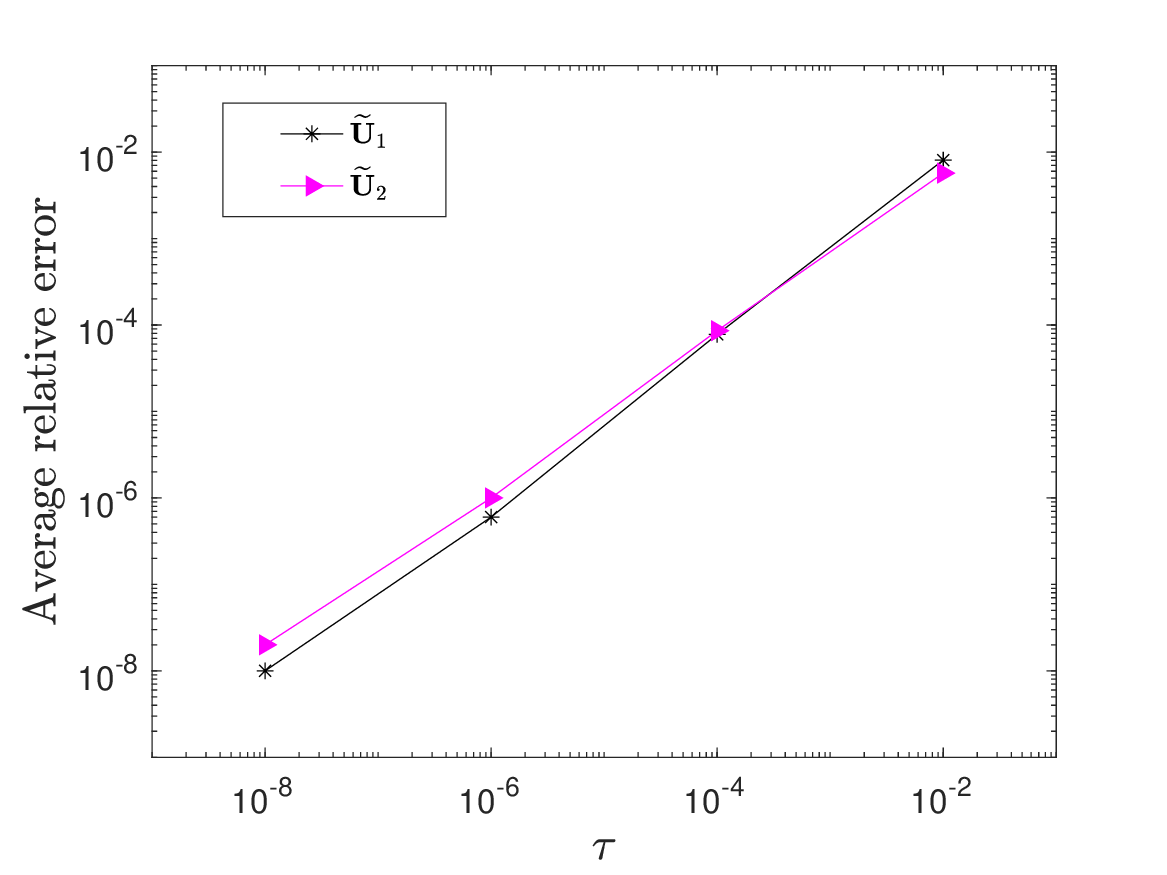}
\endminipage\hfill
\minipage{0.49\textwidth}
  \includegraphics[width=\linewidth]{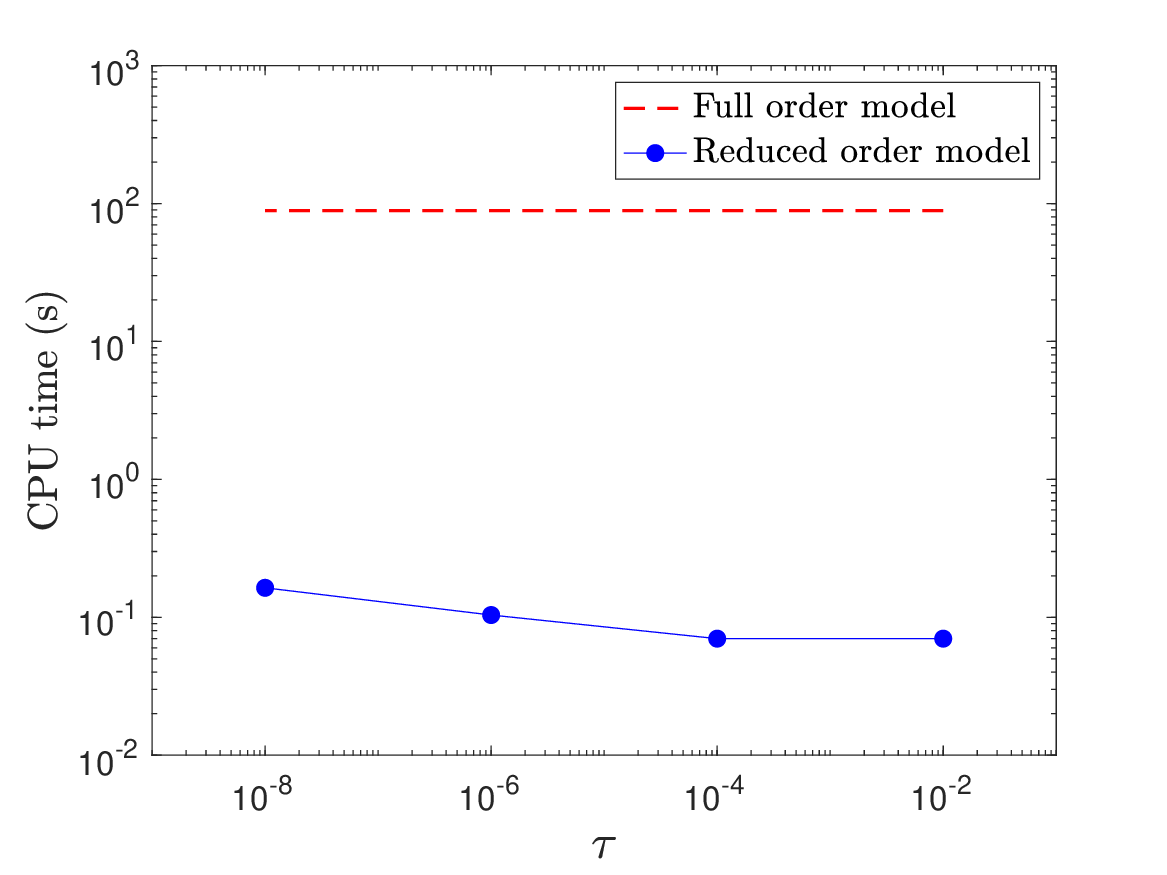}
\endminipage
    \caption{\color{black} Example \ref{ex:system}: Average relative error \ref{errormeasure2} (left) and online computational time (right) of the reduced order model and the full order model for different values of $\tau$.\label{fhnfig}}
\end{figure}
The figures indicate that even when the \MPDEIM\ reduced order model approximates the full order model with eight digits of accuracy, the time needed to integrate the model is almost three orders of magnitude faster.$\square$}

\end{example}
 In what follows we analyze the efficiency of the offline phase and compare the procedure to the standard \podeim\ procedure applied to the 2D Burgers equation in \cite{wang2016}. \begin{example}\label{ex:systemburger}
{\rm {\it The 2D {\color{black}coupled} Burgers equation (BE) \cite{wang2016}.} Here we consider the {\color{black}semilinear} 2D {\color{black}coupled} Burgers equation given by
\begin{equation}
\begin{cases}
\label{hamburger}
\dot{u}_1 =&  \frac{1}{r}\Delta u_1 - u_1 (u_1)_x - u_2 (u_1)_y\\
\dot{u}_2  =&  \frac{1}{r}\Delta u_2  - - u_1 (u_2)_x - u_2(u_2)_y
\end{cases}
\end{equation}
where $u_1(x,y,t)$ and $u_2(x,y, t)$ ($t \in [0,1]$) are the velocities to be determined, with $(x,y) \in [0,1]^2$, and $r$ is the Reynold's number. As is done in \cite{wang2016}, we derive the initial and boundary conditions from the exact traveling wave solution of the 2D Burgers equation, given by (see e.g., \cite{fletcher1983})
 \begin{equation*}
\begin{split}
u_1(x, y, t) = \frac{3}{4} - \frac{1}{4} \left(1 + e^{\frac{r(-4x + 4y -t)}{32}}   \right)^{-1} \qquad u_2(x, y, t) =  \frac{3}{4} + \frac{1}{4}\left(1 + e^{\frac{r(-4x + 4y -t)}{32}}   \right)^{-1}.
\end{split}
\end{equation*}
We consider the case $r = 100$ and discretize the model on a grid with $n$ spatial nodes in each direction, yielding a system of the form \ref{arraybigsystem}, with nonlinear functions
\begin{equation}\label{burgernonlin}
{\mathcal{F}}_i({\mathcal{D}}_{i}(\pmb{\mathcal{ U}}_i), \pmb{\mathcal{ U}}_1, \pmb{\mathcal{ U}}_2, t) = {\mathcal{F}}_i({\mathcal{D}}_{i}({{\bf U}}_i), {\bf U}_1, {\bf U}_2, t) := ({\bf B}_{1i}{\bf U}_i)\circ {\bf U}_1  +  ({\bf U}_i{\bf B}_{2i}^{\top})\circ{\bf U}_2,
\end{equation}
for $i = 1,2$, where the matrices ${\bf B}_{1i} \in \RR^{n \times n}$ and ${\bf B}_{2i} \in \RR^{n \times n}$ contain the coefficients for a first order centered difference space discretization in the $x-$ and $y-$ directions respectively ({\color{black} i.e., ${\bf B}_{1i} = {\bf B}_{2i} = \frac{n-1}{2}\, \mbox{tridiag}(-1,\underline{0}, 1)$)}, and $\circ$ is the matrix Hadamard product.  {\color{black} An upwind scheme can also be considered for ${\bf B}_{1i}$ and ${\bf B}_{2i}$, as is typically done for the coupled Burgers equation, however, in order to reproduce the results of \cite{wang2016} we consider centered finite differences.}
\begin{figure}[htb!]
\minipage{0.33\textwidth}
  \includegraphics[width=\linewidth]{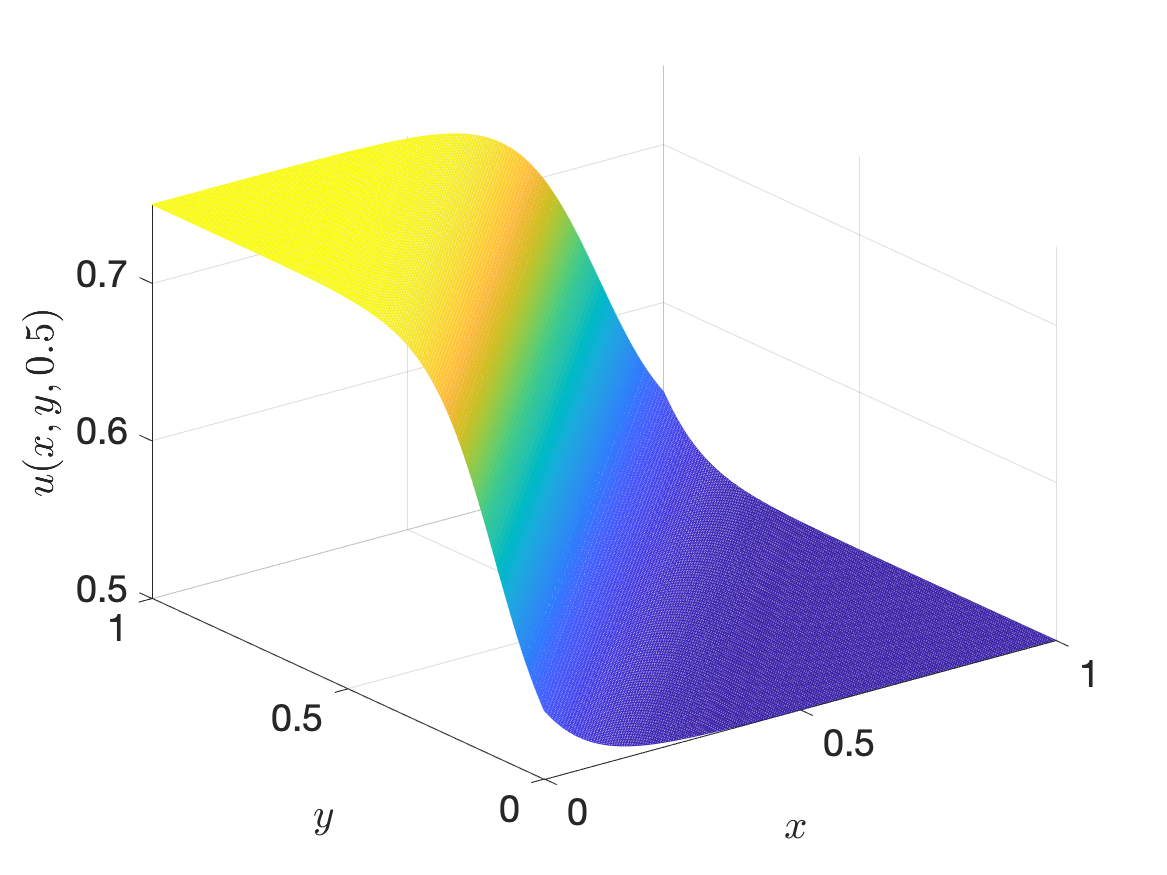}
\endminipage\hfill
\minipage{0.33\textwidth}
  \includegraphics[width=\linewidth]{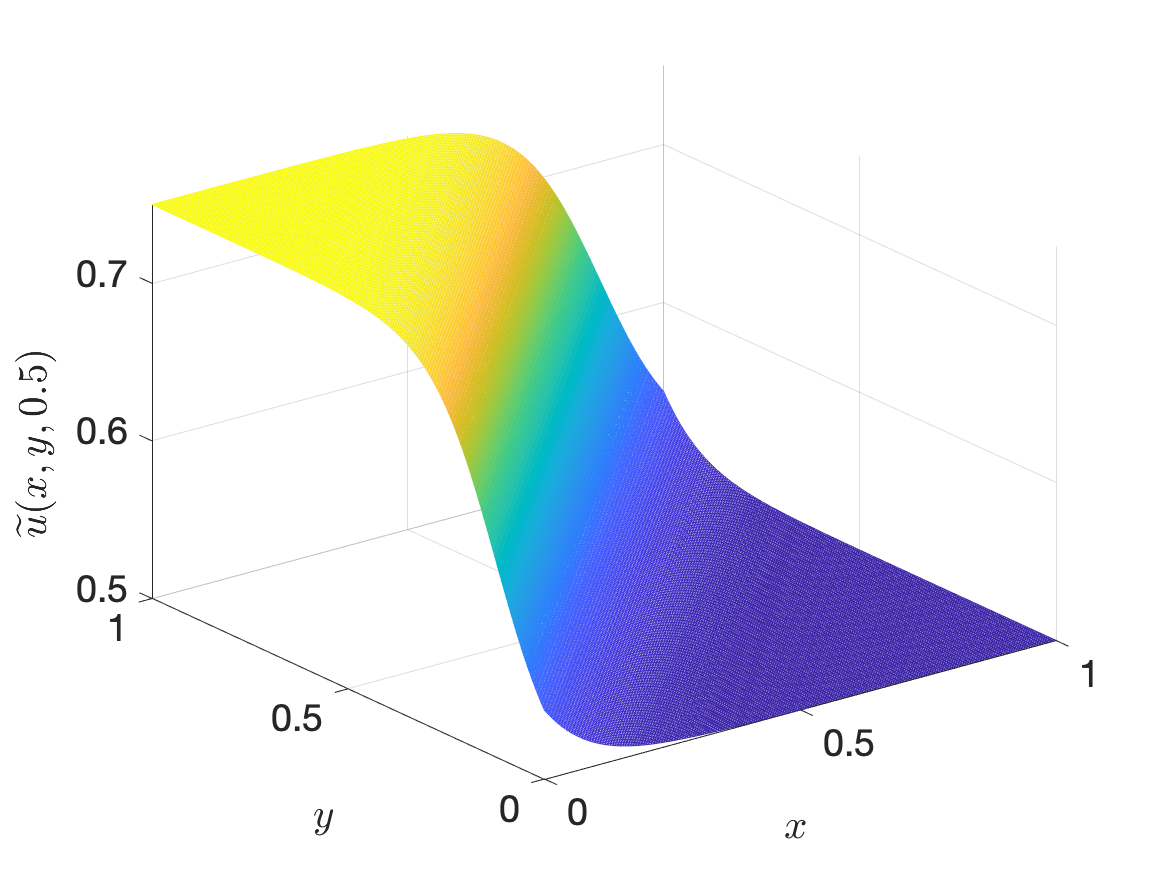}
\endminipage
\minipage{0.33\textwidth}
  \includegraphics[width=\linewidth]{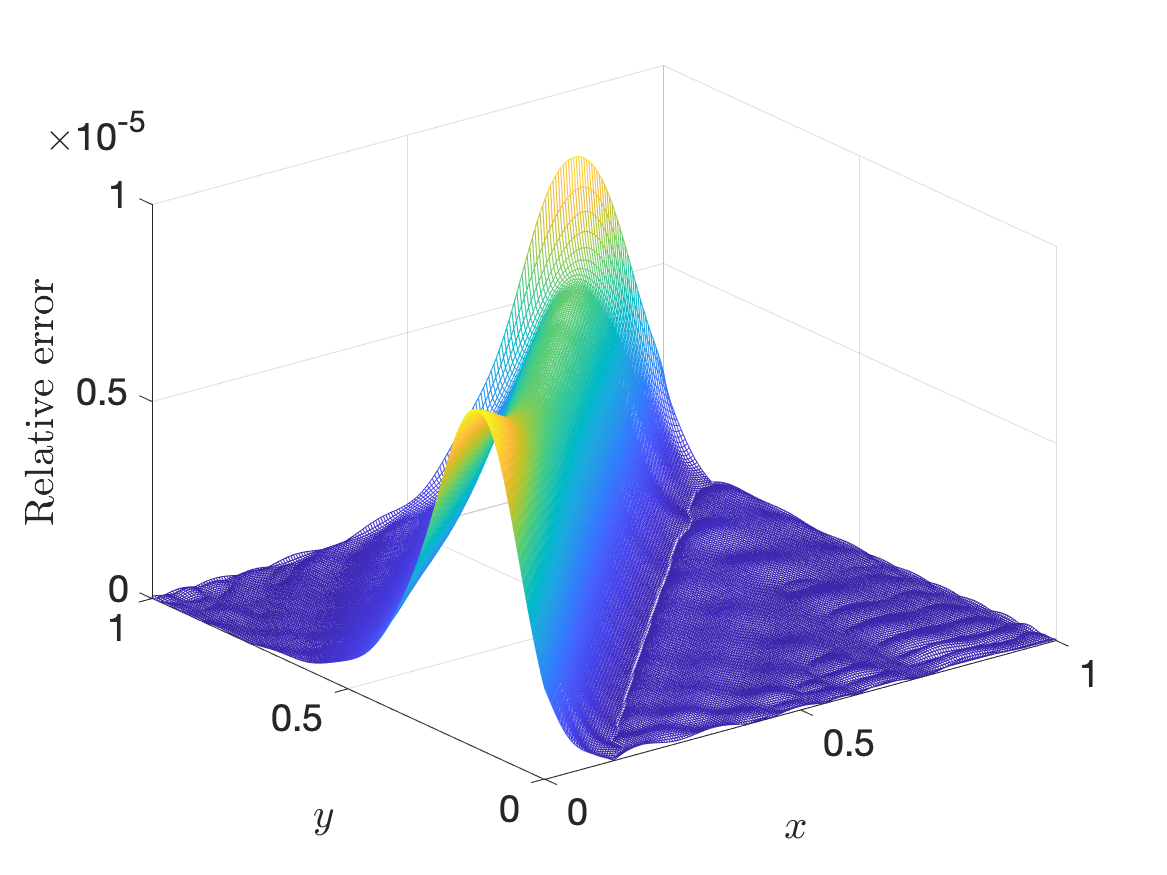}
\endminipage
    \caption{Example \ref{ex:systemburger}: $u_1(x,y,0.5)$ discretized with $n = 200$. The exact solution (left), the {\sc ho-pod-deim} approximation (middle), and the relative error mesh between the two (right).\label{ucompare}}
\end{figure}

{\color{black} It is also worth motivating the use of DEIM for the nonlinear function \ref{burgernonlin}. Indeed, this type of nonlinearity can also be efficiently treated in the vectorized POD reduced model by writing it as a tensor (see e.g., \cite{kunisch1999}) in order to avoid the use of DEIM. Nevertheless, in \cite[Table \rom{1}]{cstefuanescu2014} this idea is compared to that of \podeim\ and it is concluded that for quadratic nonlinearites the \podeim\ model requires considerably fewer floating point operations online for moderate DEIM dimension $p$. Furthermore, we use DEIM to reproduce and compare to the results of \cite{wang2016}, where it is indeed used.}

{\color{black}As mentioned above,} the presented {\sc ho-pod-deim} order reduction strategy is compared to the standard {\sc pod-deim} applied to \ref{hamburger} in \cite{wang2016}. We did not have access to the codes of \cite{wang2016}, but the \podeim\ algorithm was implemented as discussed in their paper and the results in terms of basis dimension to accuracy are comparable to the ones reported in \cite{wang2016}. Moreover, in \cite{wang2016}, the reduced order model is integrated by a fully implicit scheme, whereas for this experiment we use the IMEX 2-SBDF method to integrate both the \MPDEIM\ and \podeim\ reduced order models, which accounts for the faster online phase for {\sc pod-deim} in comparison to the times reported in \cite{wang2016}.

To this end, we consider four different space discretizations ($n = 60,200,600,1200$) and compare the computational details of \MPDEIM\ to that of \podeim\ \cite{wang2016}. Therefore, to ensure stability in the numerical integration we consider $n_{\mathfrak t} = 2n$ discrete timesteps for integrating the reduced order model. Moreover, to correspond with the space discretization error, we set $\tau = 1/n^2$.

 In our experiments, we have observed that the {\sc ho-pod-deim} strategy requires far fewer snapshots than {\sc pod-deim} to construct an equally accurate reduced order model. We hypothesize that this is because one matrix snapshot contains information about several spatial directions in $\RR^{n}$, whereas one vector snapshot only offers information about one spatial direction in $\RR^{N}$. Therefore, to obtain a {\sc pod} basis of dimension $k$ in the vectorized setting, at least $k$ snapshots are required, even though one vector in $\RR^{N}$ contains many spatial directions from $\RR^n$. A further experimental or theoretical analysis of the relationship between array snapshots and vectorized snapshots could be interesting in future work. To this end, we consider $n_s = 20$ equispaced snapshots for {\sc ho-pod-deim} and $n_s = 100$ for {\sc pod-deim}. 
\begin{figure}[htb]
\minipage{0.33\textwidth}
  \includegraphics[width=\linewidth]{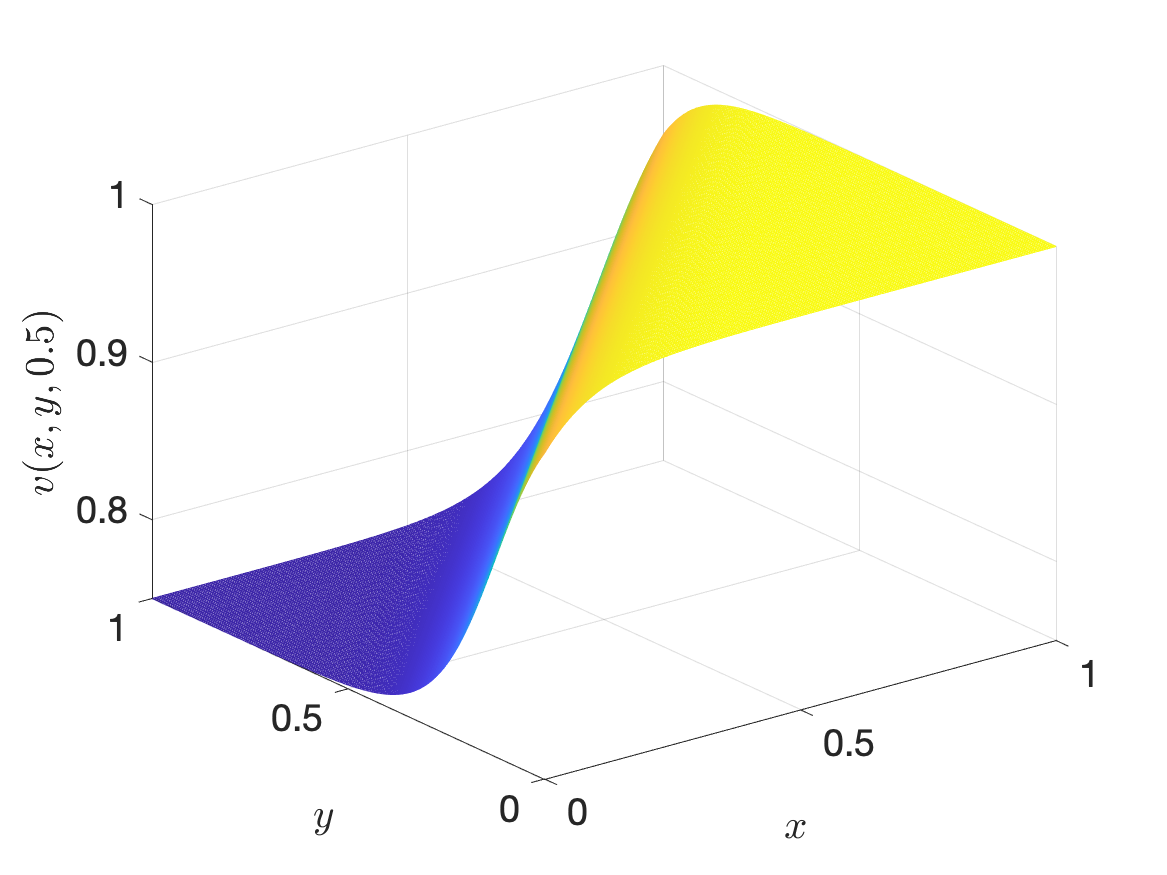}
\endminipage\hfill
\minipage{0.33\textwidth}
  \includegraphics[width=\linewidth]{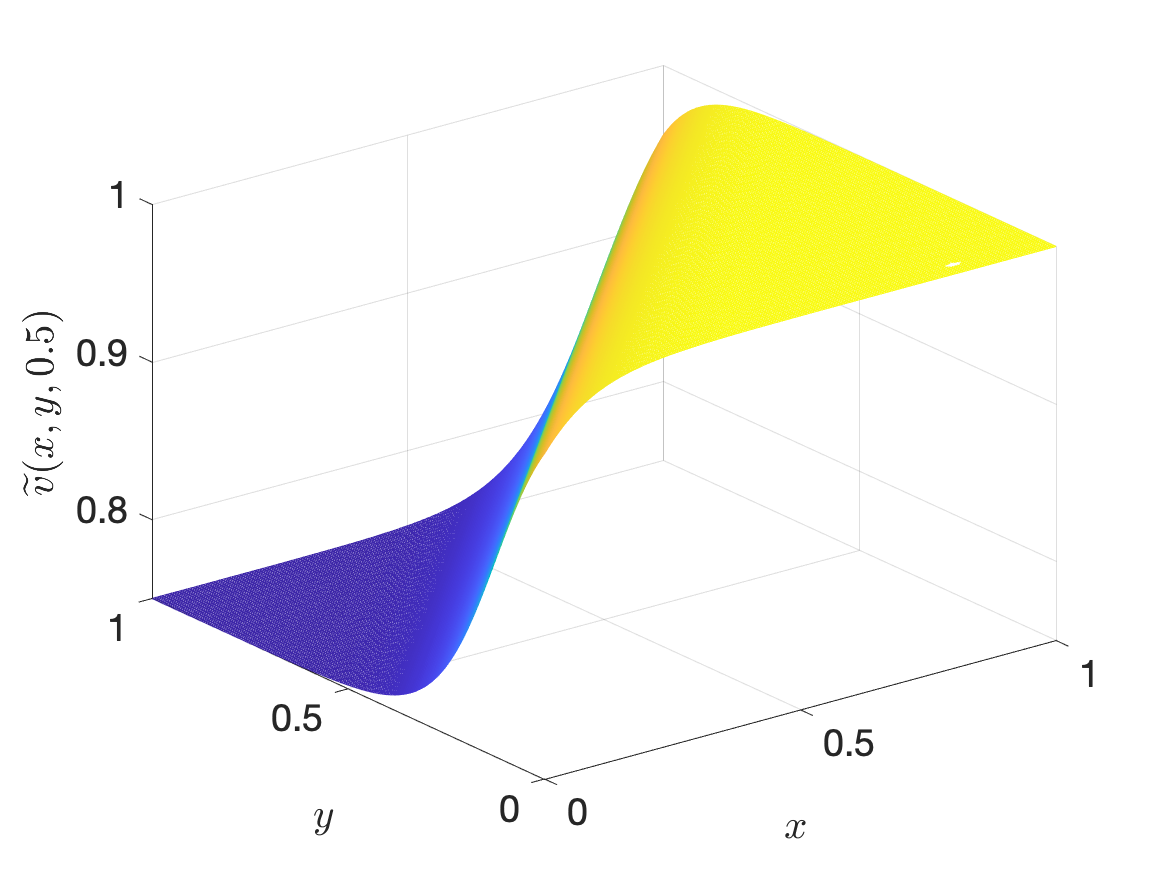}
\endminipage
\minipage{0.33\textwidth}
  \includegraphics[width=\linewidth]{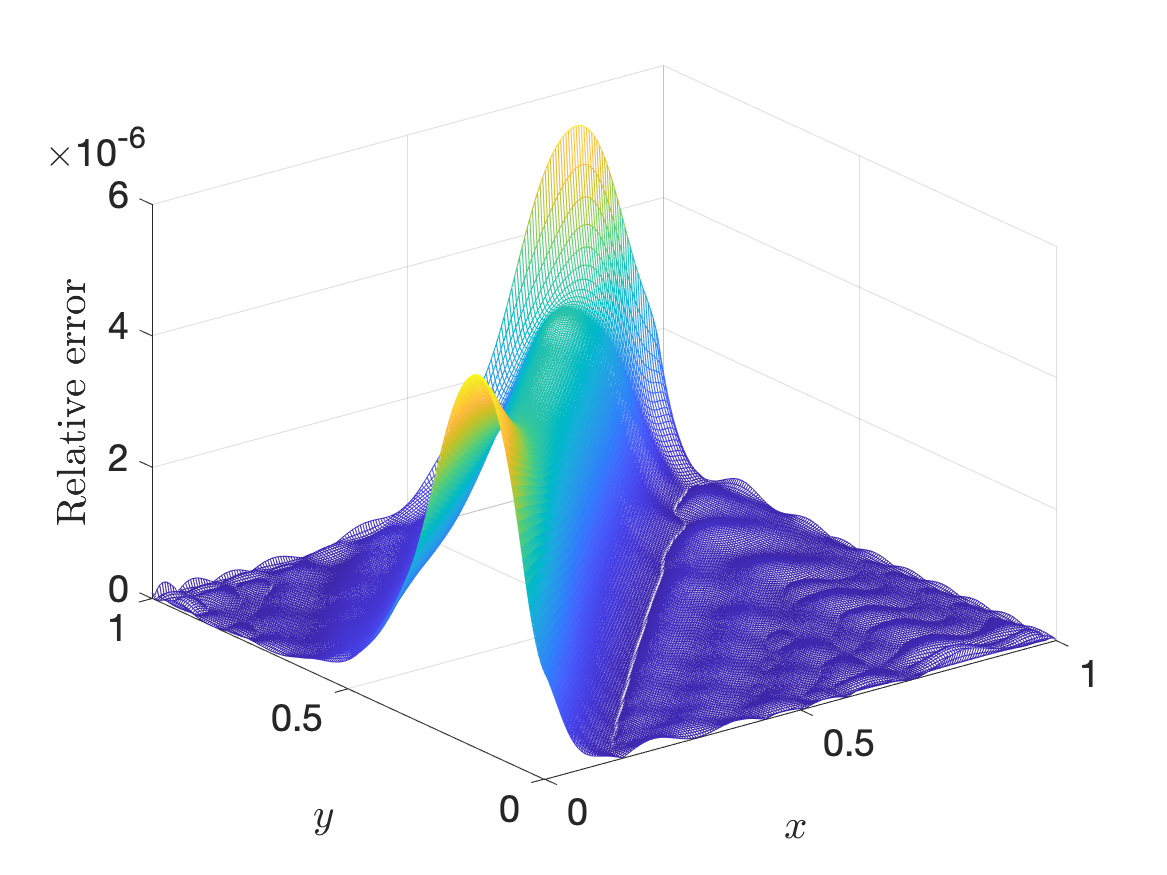} 
\endminipage
    \caption{Example \ref{ex:systemburger}:  $u_2(x,y,0.5)$ discretized with $n = 200$. The exact solution (left), the {\sc ho-pod-deim} approximation (middle), and the relative error mesh between the two (right).\label{vcompare}}\end{figure}

A visual comparison of the accuracy of the {\sc ho-pod-deim} reduced model at $t = 0.5$, when $n = 200$ is plotted in Figures \ref{ucompare} and \ref{vcompare} for $u_1$ and $u_2$ respectively. Moreover, in Figure \ref{u1errors} we plot the average relative error through $n_{\mathfrak t} = 2n$ timesteps between the {\sc ho-pod-deim} approximation and the exact solution at the relevant nodes, for the four different space dimensions $n$. We investigate the computational load required by both strategies to achieve this accuracy.

 Firstly, we report in Table \ref{dimtable} the reduced basis dimensions for {\sc ho-pod-deim} and {\sc pod-deim} for all four space discretizations, as well as the memory requirements. In particular, for each ${\bf U}_i$ we report the dimensions $k_1/k_2$ ($p_1/p_2$) of the {\sc ho-pod} ({\sc ho-deim}) bases and the dimension $k$ ($p$) of the {\sc pod} ({\sc deim}) bases. Moreover, the reported global memory requirements include the number of stored vectors 
multiplied by their length ($\#\cdot {length}$) in each phase. 
\begin{figure}[htb]
\centering
  \includegraphics[width=0.5 \linewidth]{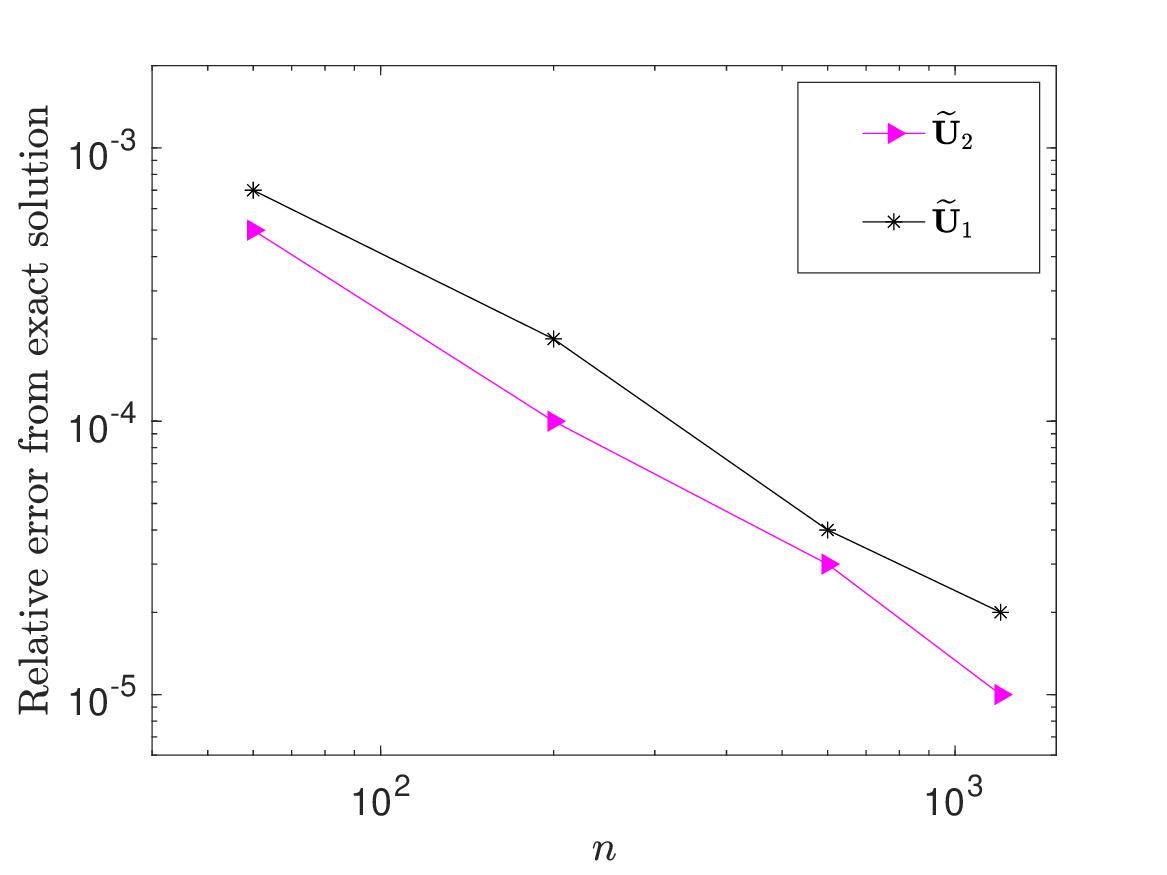}
    \caption{Example \ref{ex:systemburger}: The average relative error through $n_{\mathfrak t} = 2n$ timesteps between the \MPDEIM\ approximation $\widetilde{\bf U}_1(t)$ ($\widetilde{\bf U}_2(t)$) and the exact solution $u_1(x,y,t)$ ($u_2(x,y,t)$).\label{u1errors}}
\end{figure}
The table indicates a great reduction in memory requirements for {\sc ho-pod-deim} in comparison to {\sc pod-deim}, but the {\sc pod-deim} strategy produces a smaller reduced order model, as is expected from a one-sided reduction strategy. That is, the {\sc pod-deim} reduced model requires evaluating the nonlinear function at merely $p$ entries, whereas the \MPDEIM\ reduced model requires evaluating the nonlinear function at $p_1\cdot p_2$ entries. From Table \ref{dimtable} it is clear that $p_1 \cdot p_2 \ge p$ for all $n$.  

We investigate the pros and cons, in terms of computational time, of both strategies in Figure \ref{offandon}.
 \begin{table}[]
\centering
\caption{A breakdown of the {\sc (ho)-pod} and {\sc (ho)-deim} basis dimensions and the memory requirements for four different state space dimensions. Note that $\tau = 1/n^2$. \label{dimtable}}
\begin{tabular}{|l|r|r|r|r|c|c|}
\hline
          &&   &   & & \sc offline &\sc online  \\ 
$n$&  {\sc algorithm}  &  ${\bf U}_i$&\sc pod dim. & \sc deim dim. & \sc memory & \sc memory \\
\hline  
{\multirow{4}{*}{$60$}} & \multirow{2}{*}{{\sc ho-pod-deim}}   & ${\bf U}_1$ & 9/9                                  & 18/18                                 & $98 n$                          & $54n$                         \\ 
\multicolumn{1}{|c|}{}                          &                                                       & ${\bf U}_2$ & 9/9                                  & 18/18                                 & $98n$                          & $54n$                         \\ \cline{2-7} 
\multicolumn{1}{|c|}{}                          & \multirow{2}{*}{{\sc pod-deim} \cite{wang2016} }    & ${\bf U}_1$ & 5                                    & 14                                    & $400n^2$                        & $19n^2$                       \\ 
\multicolumn{1}{|c|}{}                          &                                                       & ${\bf U}_2$ & 4                                    & 14                                    & $400n^2$                        & $18n^2$                       \\ \hline 
\multirow{4}{*}{$200$}                      & \multirow{2}{*}{{\sc ho-pod-deim}}   & ${\bf U}_1$ & 13/13                                & 24/25                                 & $153n$                          & $75n$                         \\
                                                &                                                       & ${\bf U}_2$ & 12/12                                & 24/25                                 & $153n$                          & $73n$                         \\ \cline{2-7} 
                                                & \multirow{2}{*}{{\sc pod-deim} \cite{wang2016}}      & ${\bf U}_1$ & 9                                    & 23                                    & $400n^2$                        & $32n^2$                       \\
                                                &                                                       & ${\bf U}_2$ & 8                                    & 23                                    & $400n^2$                        & $31n^2$                        \\ \hline 
\multirow{4}{*}{$600$}                      & \multirow{2}{*}{{\sc ho-pod-deim}}   & ${\bf U}_1$ & 16/17                                & 32/32                                 & $196n$                          & $97n$                         \\
                                                &                                                       & ${\bf U}_2$ & 16/16                                & 32/32                                 & $194n$                          & $96n$                         \\ \cline{2-7} 
                                                & \multirow{2}{*}{{\sc pod-deim}  \cite{wang2016} }     & ${\bf U}_1$ & 15                                   & 28                                    & $400 n^2$                        & $43n^2$                       \\ 
                                                &                                                       & ${\bf U}_2$ & 14                                   & 28                                    & $400n^2$                        & $42n^2$                        \\ \hline 
\multirow{4}{*}{$1200$}                                      & \multirow{2}{*}{{\sc ho-pod-deim}}   & ${\bf U}_1$ & 19/19                                & 36/39                                 & $219n$                          & $113n$                        \\ 
                                                &                                                       & ${\bf U}_2$ & 19/19                                & 36/39                                 & $215n$                          & $113 n$                        \\\cline{2-7} 
                                                &  \multirow{2}{*}{{\sc pod-deim}  \cite{wang2016} }                  & ${\bf U}_1$ & 19                                   & 31                                    & $400n^2$                        & $50 n^2$                       \\
                                                &                                                       & ${\bf U}_2$ & 18                                   & 31                                    & $400 n^2$                        & $50 n^2$                       \\ \hline
\end{tabular}
\end{table}
On the left of Figure~\ref{offandon} we plot the time needed offline to construct the basis vectors for both strategies, for increasing $n$. For {\sc ho-pod-deim} this includes the time needed to perform the SVD of each snapshot and the time needed to orthogonalize and truncate the new basis vectors for all 8 bases\footnote{Each equation $u_1$ and $u_2$ require four basis matrices when $d = 2$. Two stemming from the snapshot solutions $\{ \pmb{\mathcal{ U}}_i (t_{j})\}_{j = 1}^{n_s}$ for the \MPOD\ dimension reduction and two stemming from the nonlinear snapshots \ref{nonsnaps} for the \MODEIM\ interpolation.}, whereas for {\sc pod-deim} this includes the time to vectorize each snapshot and the time to perform the economy SVD of all four  $n^2 \times n_s$ matrices of snapshots. On the right of Figure~\ref{offandon} we report the time needed to evaluate \ref{arraysmallsystem} and a system of the form \ref{redode} at $n_{\mathfrak t}$ timesteps online and compare it to the time needed to evaluate the full order model \ref{arraybigsystem}. The timings in Figure~\ref{offandon} (right) correspond to the reduced dimensions reported in Table \ref{dimtable}.

Figure~\ref{offandon} (left) indicates the large gain in offline computational time by the new strategy, with almost two orders of magnitude difference as $n$ increases. However, due to the larger reduced dimensions of {\sc ho-pod-deim} reported in Table \ref{dimtable}, fractionally more time is required online, as presented in Figure~\ref{offandon} (right). Nevertheless the online times are very comparable and orders of magnitude lower than the full order model \ref{arraybigsystem}. In problems where the nonlinear term is more expensive to evaluate this drawback of {\sc ho-pod-deim} will become more evident, however, and will need to be further investigated in future work; {\color{black} see also Remark~\ref{deimremark}}.
%
\begin{figure}[htb]
\minipage{0.49\textwidth}
  \includegraphics[width=\linewidth]{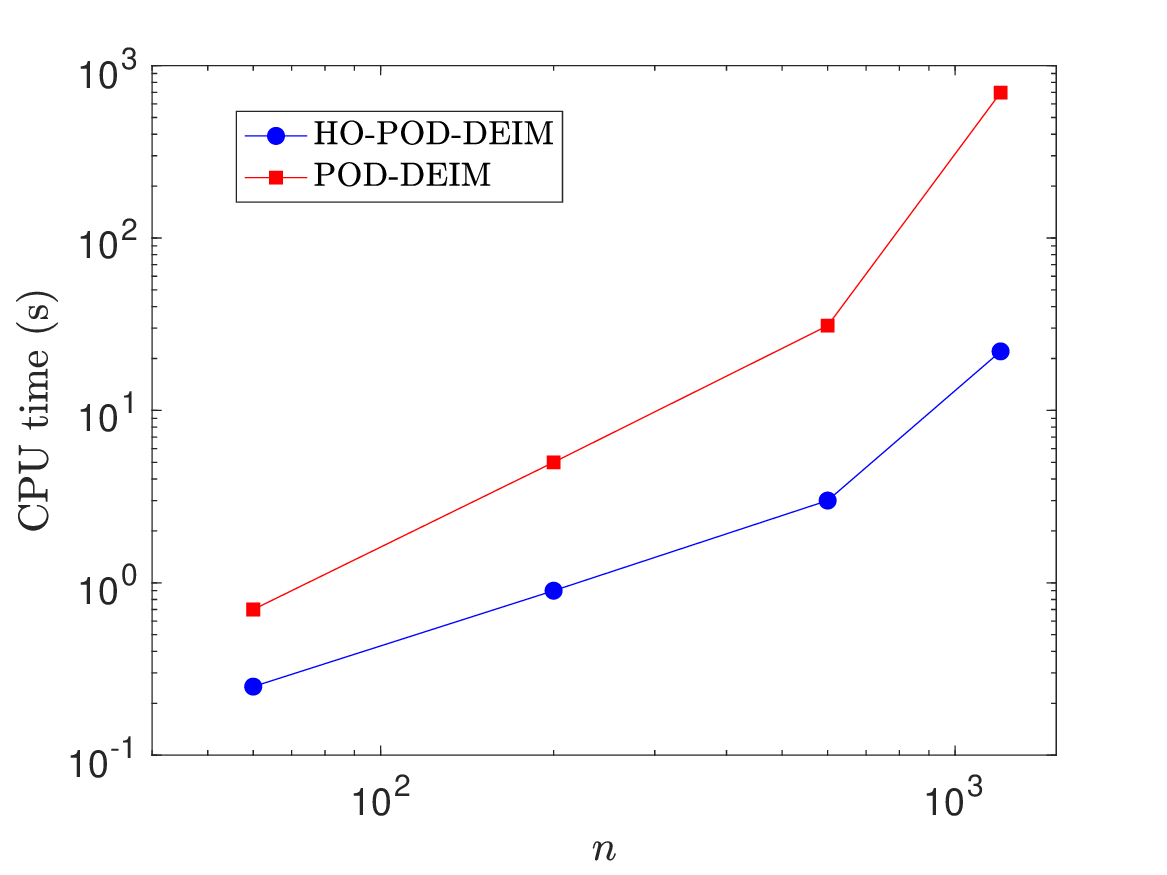}
\endminipage\hfill
\minipage{0.49\textwidth}
  \includegraphics[width=\linewidth]{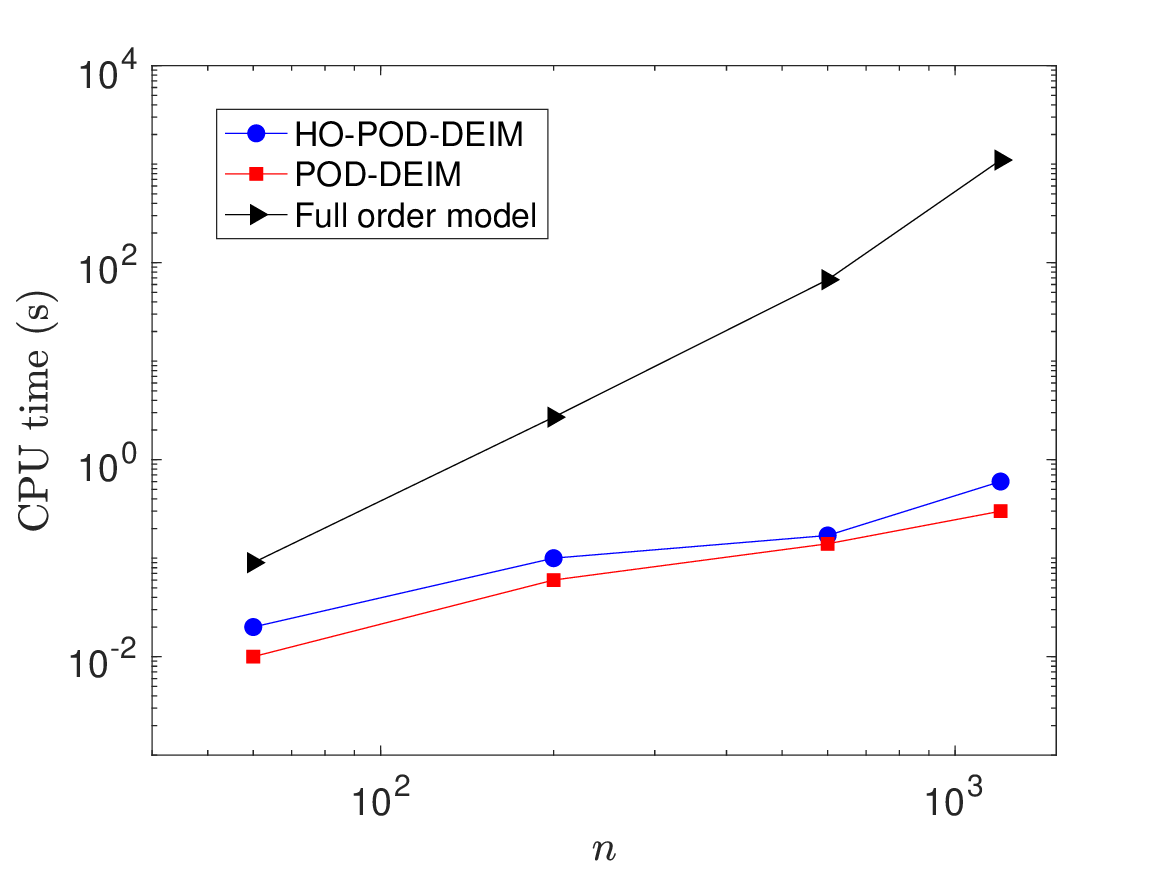}
\endminipage
    \caption{Example~\ref{ex:systemburger}: A comparison of the time required offline for basis construction (left) and online for integration (right) between {\sc ho-pod-deim} and {\sc pod-deim} \cite{wang2016} for increasing $n$.\label{offandon}}
\end{figure}

This experiment indicated that a greater accuracy with respect to the exact solution can be achieved by the discrete {\sc ho-pod-deim} reduced order model in a fraction of the offline computational time compared to {\sc pod-deim}. Moreover the online time remains comparable, and a large gain in memory requirements is witnessed. 
}$\square$
\end{example}
In what follows we illustrate the efficiency of the procedure in the multilinear setting.

\begin{example}\label{ex:systemburger3d}
{\rm {\it The 3D {\color{black}coupled} Burgers equation (BE).} Here we consider the {\color{black}semilinear} 3D {\color{black}coupled} Burgers equation (see, e.g., \cite{GAO2017}) given by
\begin{equation}
\begin{cases}
\label{hamburger3d}
\dot{u}_1 &=  \frac{1}{r}\Delta u_1 -  \underline{u} \cdot \nabla u_1\\
\dot{u}_2 &=  \frac{1}{r}\Delta u_2 - \underline{u} \cdot \nabla u_2,\\
\dot{u}_3 &=  \frac{1}{r}\Delta u_3 - \underline{u} \cdot \nabla u_3,\\
\end{cases}
\end{equation}
where $u_1(x,y,z,t)$, $u_2(x,y,z,t)$ and $u_3(x,y,z,t)$ are the three velocities to be determined, with $(x,y,z) \in [0,1]^3$ and $t \in [0,1]$. Furthermore, the system is subject to homogeneous Dirichlet boundary conditions and initial states
\begin{equation*}
\begin{split}
u_1(x, y, z, 0) &= \frac{1}{10}\sin(2\pi x)\sin(2\pi y)\cos(2\pi z) \\
u_2(x, y, z, 0) &=  \frac{1}{10}\sin(2\pi x)\cos(2\pi y)\sin(2\pi z)\\
u_3(x, y, z, 0) &=  \frac{1}{10}\cos(2\pi x)\sin(2\pi y)\sin(2\pi z). \\
\end{split}
\end{equation*}
A finite difference space discretization inside the cube yields a system of ODEs of the form \ref{arraybigsystem}, with nonlinear functions given by
\begin{equation*}
\begin{split}
{\mathcal{F}}_i({\mathcal{D}}_i(\pmb{\mathcal{ U}}_i),\pmb{\mathcal{ U}}_1, \pmb{\mathcal{ U}}_2,\pmb{\mathcal{ U}}_3, t) &= (\pmb{\mathcal{ U}}_i \times_1 {\bf B}_{1i}) \circ \pmb{\mathcal{ U}}_1 
+(\pmb{\mathcal{ U}}_i \times_2 {\bf B}_{2i}) \circ \pmb{\mathcal{ U}}_2
+(\pmb{\mathcal{ U}}_i \times_3 {\bf B}_{3i}) \circ \pmb{\mathcal{ U}}_3,
\end{split}
\end{equation*}
for $i = 1,2,3$, where ${\bf B}_{1i} \in \RR^{n \times n}$, ${\bf B}_{2i} \in \RR^{n \times n}$ and ${\bf B}_{3i} \in \RR^{n \times n}$ contain the coefficients for a first order centered difference space discretization in the $x-$, $y-$ and $z-$ directions respectively. Furthermore we vary $r$ through the experiment, and calculate the reduced order model through $n_s = 50$ equispaced snapshots for both \MPDEIM\ and \podeim.

Firstly, we set $r = 10$ and consider five different space discretizations ($n = 50, 80, 100, 150, 200$) and investigate the efficiency of the offline phase of the {\sc ho-pod-deim} procedure for systems when $d = 3$, in comparison to standard {\sc pod-deim}. Note that the system \ref{arraybigsystem} of dimension $n = 200$ is equivalent to the system \ref{systemode} with dimension $N = 8\, 000\, 000$ when $d=3$.

The improvement in memory requirements is immediately evident, since the new procedure requires storing basis vectors of length $n$, whereas the vectorization procedure needs to store many basis vectors of length $n^3$, as has been witnessed in Table \ref{dimtable} for $d = 2$.
\begin{figure}[htb]
\minipage{0.49\textwidth}
  \includegraphics[width=\linewidth]{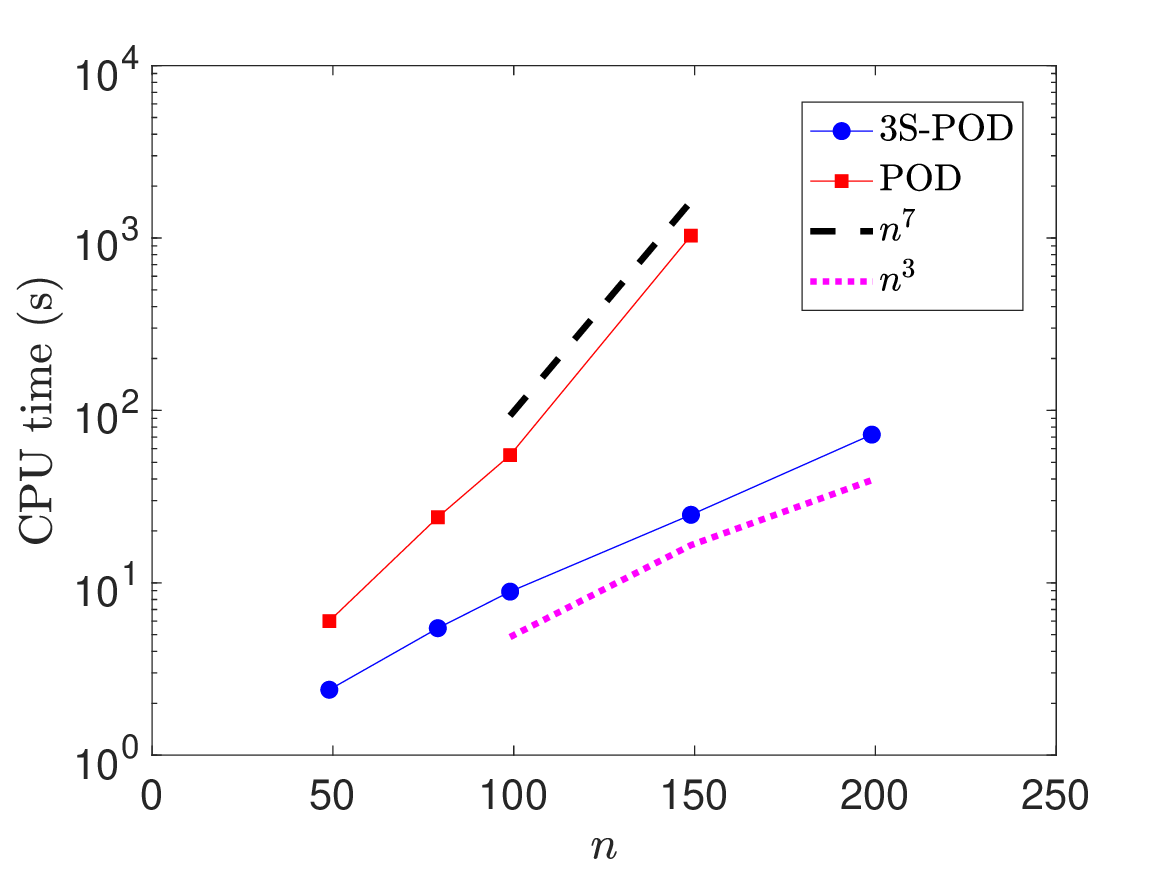}
\endminipage\hfill
\minipage{0.49\textwidth}
  \includegraphics[width=\linewidth]{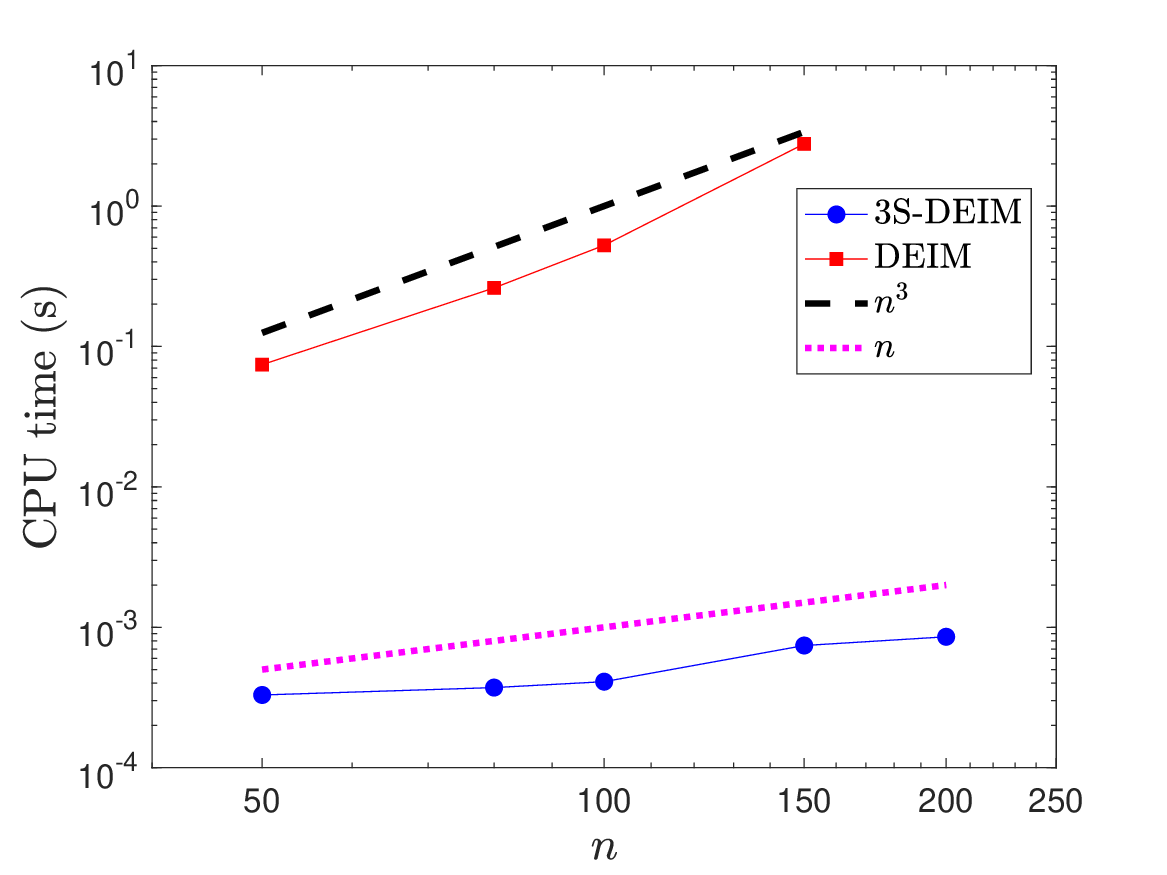}
\endminipage
    \caption{Example~\ref{ex:systemburger3d}: A comparison of the offline time for increasing dimension $n$, between \MPOD\ and {\sc pod} (left) {\color{black} and \MODEIM\ and {\sc deim} (right)}.\label{burgertimes3d}}
\end{figure}
In Figure~\ref{burgertimes3d} (left) we compare the computational time needed to determine the basis vectors, given $\tau = 10^{-4}$, for increasing dimension $n$. For the new procedure that includes the time needed to perform the STHOSVD of each snapshot and the time needed to orthogonalize and truncate the new basis vectors, whereas for {\sc pod-deim} this includes the time to vectorize each snapshot and the time to perform the economy SVD of the $n^3 \times n_s$ matrix of snapshots. The times are added together for all 18 bases\footnote{Each equation $u_1$,$u_2$ and $u_3$ requires six basis matrices when $d = 3$. Three stemming from the snapshot solutions $\{ \pmb{\mathcal{ U}}_i (t_{j})\}_{j = 1}^{n_s}$ for the \MPOD\ dimension reduction and three stemming from the nonlinear snapshots \ref{nonsnaps} for the \MODEIM\ interpolation. For standard {\sc pod-deim} each equation requires only two basis matrices, one for dimension reduction and one for {\sc deim} interpolation, hence six bases in total.} required by {\sc ho-pod-deim} and all 6 bases required by {\sc pod-deim}. {\color{black} We explicitly remark that both the STHOSVD and the economy SVD can potentially be further accelerated by using randomized algorithms; see e.g., \cite{halko2011, minster2020}. This is, however, not considered in our experiments.} 

In Figure~\ref{burgertimes3d} (right) we compare the time needed to determine the {\sc (ho)-deim} interpolation indices. That is, the cumulative time taken by {\tt q-deim} for all 9 nonlinear bases for  {\sc ho-pod-deim} and all 3 nonlinear bases required by {\sc pod-deim}.

For both {\sc pod} and {\sc deim} the improvement in computational time is very evident in the plots, with a few orders of magnitude difference. On the standard laptop computer on which these experiments were performed, the \MPDEIM\ bases were created in just more than a minute for $n = 200$, whereas for the vectorization procedure, the computer ran out of memory, after processing 30 snapshots in more than an hour.

In what follows we investigate the online phase. We set $n = 150$ and illustrate the efficiency of the three-sided reduction procedure, together with the new {\sc t3-sylv} method for solving the low-dimensional, dense tensor-valued system of equations. To this end we investigate the total time needed for solving all inner linear systems at $n_{\mathfrak t} = 100$ timesteps, that is $300$ linear systems in total, using {\sc t3-sylv} and compare it to the time needed if the system \ref{arraysystemtime} is vectorized and solved as a standard (nearly dense) linear system (Vec-lin). This is done for different values of $\tau$, which result in different reduced dimensions. In particular, we plot the computational time with respect to the maximum dimension of the vectorized systems for the different values of $\tau$. That is, the value $\mbox{max}(k_1k_2k_3)$ on the $x-$axis is the maximum value value of $k_{1i}k_{2i}k_{3i}$ for all $i = 1,2,3$, given $\tau$.

For the solution of the vectorized system we perform a  reverse Cuthill–McKee reordering of the coefficient matrices, to exploit any remaining sparsity pattern, and perform an $LU$ decomposition once for all, so that only front- and back-substitution is required for all system solves. The results are reported in Figure~\ref{innertimes}.
\begin{figure}[htb]
\centering
  \includegraphics[width=0.5\textwidth]{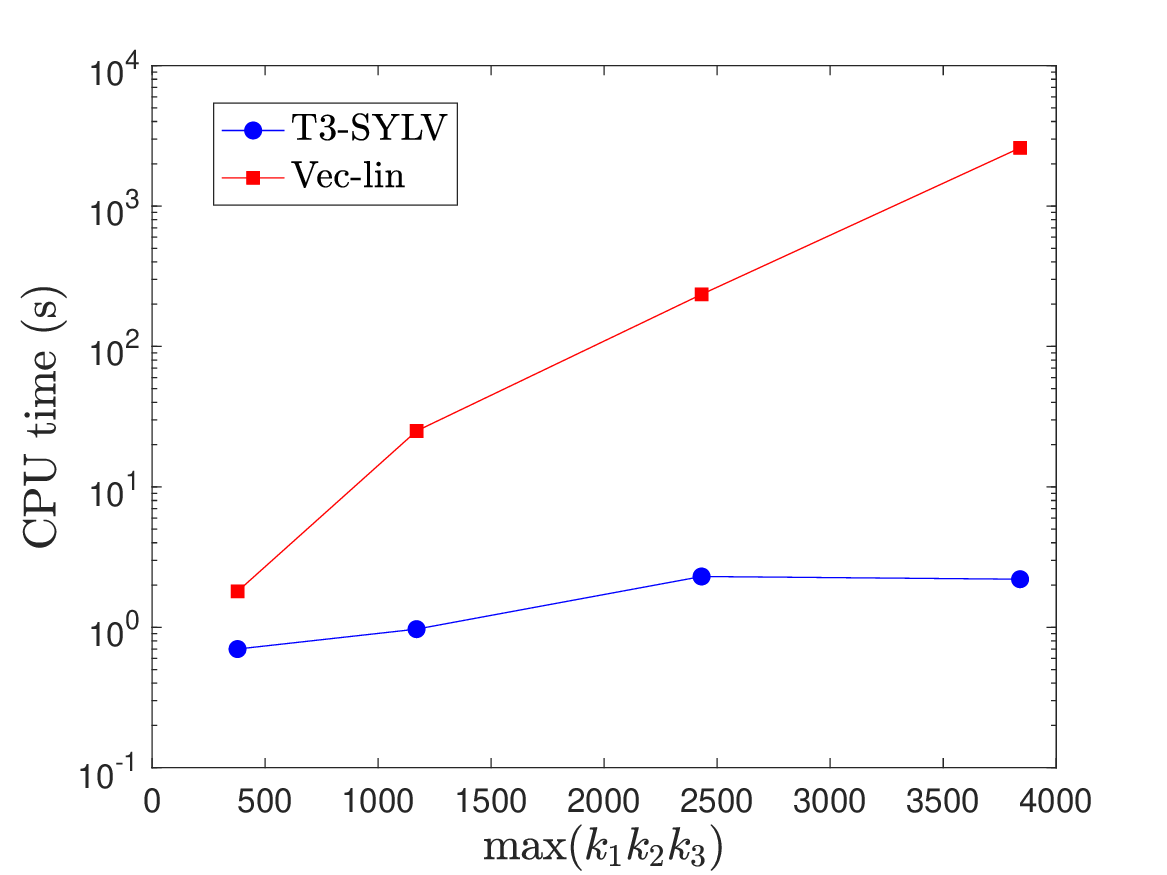}
    \caption{Example~\ref{ex:systemburger3d}: A comparison of the time to solve all linear systems of the form \ref{arraysystemtime}, for different values of $\tau$, between {\sc t3-sylv} and Vec-lin. The $x-$axis displays the maximum dimension of the three vectorized equations for different values of $\tau$. \label{innertimes}}
\end{figure}

The advantage that the three-sided reduction procedure poses in combination with the {\sc t3-sylv} inner solver is evident from Figure~\ref{innertimes}. This figure, together with Figure~\ref{burgertimes3d} illustrates that a lot of time can be saved offline and online. Without the {\sc t3-sylv} inner solver it is evident that the tensor structure of the coefficient matrices retained by the three sided projection would result in expensive, dense linear system solves, which would cancel the time that has been saved in the offline phase. Nevertheless, with the availability of the {\sc t3-sylv} solver, a large speedup is seen in both the offline and online phases, together with a massive gain in memory requirements.

Finally, Table~\ref{comp5} and Table~\ref{comp6} contain the details of the reduced order model, given $n = 150$, $\tau = 10^{-4}$ and the error measure \ref{errormeasure2} for increasing values of the Reynold's number $r$.
\begin{table}[bht]
\caption{Example~\ref{ex:systemburger3d}. Dim. of \MPOD\ and \MODEIM\ bases and  the average error at 300 timesteps for increasing $r$. The full order model has dimension $n = 150$ and $\tau = 10^{-4}$.  \label{comp5}}
\centering
\begin{tabular}{|l|l|c|c|c|c|c|c|c|}
\hline             
     &    &        &    &     &       &   &  & {\sc error} \\
    $r$ &   $u$  &   $k_1$       & $k_2$       &$k_3$    &   $p_1$        &  $p_2$     &$p_3$      & $\bar{\mathcal{E}}(\pmb{\mathcal{ U}})$ \\ \hline
   \multirow{3}{*}{{$10$}} &$u_1$  & 4        & 7   & 10        & 7& 12& 16        & $1\cdot 10^{-4}$   \\ 
&$u_2$  & 7          & 7  & 7        & 9 &12&13        & $6\cdot 10^{-5}$    \\ 
&$u_3$  & 8          & 12  & 8        & 9 &16&13        & $1\cdot 10^{-4}$    \\ \hline
   \multirow{3}{*}{{$100$}} &$u_1$  & 6        & 11   & 15        & 10 & 17& 20        & $3\cdot 10^{-5}$    \\ 
&$u_2$  & 10          & 11  & 11        & 12 &17&17        & $4\cdot 10^{-5}$    \\ 
&$u_3$  & 10         & 16  & 12        & 12 &21&17        & $4\cdot 10^{-5}$    \\ \hline
   \multirow{3}{*}{{$500$}} &$u_1$  & 9       & 15   & 19        & 13 & 23& 26        & $2\cdot 10^{-5}$   \\ 
&$u_2$  & 11          & 16  & 17        & 14 &23&23        & $3\cdot 10^{-5}$    \\ 
&$u_3$  & 12         & 19  & 16        & 14 &25&23        & $4\cdot 10^{-5}$    \\ \hline
\end{tabular}
\end{table}
Table~\ref{comp5} illustrates that a large reduction in all dimensions is achieved, with a very acceptable accuracy over 300 timesteps, even for large $r$. Nevertheless it is clear that for larger Reynold's number the singular value decay in each mode becomes slower, resulting in larger reduced dimensions. Moreover we mention that for this same problem, when $r = 100$, standard {\sc pod-deim} would require storing 57 vectors of length $n^3$ in the online phase, as opposed to the 245 vectors of length $n$ that need to be stored for \MPDEIM. Furthermore, Table~\ref{comp6} confirms that even for large $r$ we observe a large gain in online computational time achieved by the reduced model.$\square$

\begin{table}[bht]
\caption{Example~\ref{ex:systemburger3d}. Memory and CPU time required for basis construction and integration. The full order model has dimension $n = 150$  and $\tau = 10^{-4}$.  \label{comp6}}
\centering
\begin{tabular}{|l|c|c|c|c|}
\hline             
   &      Online    & Basis   & FOM & ROM   \\
  r&       memory       &time(s)     & time(s)    & time(s)  \\ \hline
    10  &  $177n$   & 20       & 1641         & 1.9  \\ \hline
    100  &  $245n$   & 20       & 1641         & 2.2 \\ \hline
     500 &   $318n$   & 20       & 1641         & 3.3  \\ \hline
\end{tabular}
\end{table}}
\end{example}
{\color{black}
\begin{example}\label{ex:cells3d}
{\rm {\it A 3D reaction-diffusion model for cell apoptosis.} 
As a final example we consider a reaction-diffusion system, orginally introduced in \cite{daub2012} to investigate the behavior of protein concentrations (in space and time) of a cell apoptosis model in 1D. The model was later extended to higher dimension in \cite[Chapter 2.3]{daub2013}. The proteins build a network called ``caspase--cascade'' and the dynamics, with homogeneous Neumann boundary conditions, are given by
{\small\begin{equation}
\begin{aligned}
\label{cell3d}
\dot{u}_1 &=  \delta_1\Delta u_1 - c_4u_1 + c_1\sin(u_3u_2), &\quad \dot{u}_2 &=   \delta_2\Delta u_2 - c_4u_2 + c_2u_4u_1^3,  \\
\dot{u}_3 &=  \delta_3\Delta u_3 - c_4u_3 - c_1\sin(u_3u_2) + c_3,&\quad \dot{u}_4 &= \delta_4\Delta u_4 - c_4u_4 - c_2u_4u_1^3 + c_3,
\end{aligned}
\end{equation}}where $u_1({\bf x},t)$, $u_2({\bf x},t)$, $u_3({\bf x},t)$ and $u_4({\bf x},t)$ are four different reactants called Procaspase-8, Procaspase-3, Caspase-8 and Caspase-3 respectively, with ${\bf x} = (x,y,z) \in [0,1]^3 =: \Omega$ and $t \in [0,1]$. For the values and derivation of the constants $\delta_i$ and $c_i$ we refer the reader to \cite[Chapter 2.3]{daub2013} and we consider the initial condition
$$
\left(u_1, u_2, u_3, u_4\right)({\bf x},0) = \begin{cases} \left(u_1^{(d)}, u_2^{(d)}, u_3^{(d)}, u_4^{(d)} \right) \quad \mbox{for} \quad &{\bf x} \in \Omega_{\tiny \mbox{ext}}\\     \left(u_1^{(\ell)}, u_2^{(\ell)}, u_3^{(\ell)}, u_4^{(\ell)} \right) \quad \,\,\mbox{for} \quad &{\bf x} \in \Omega_{\tiny \mbox{in}} \end{cases},
$$
where $\Omega_{\tiny \mbox{ext}} := \{{\bf x} \in \Omega, r_0 \le \|{\bf x}\|_2 \le 1 \}$ and $\Omega_{\tiny \mbox{in}} := \{{\bf x} \in \Omega, \|{\bf x}\|_2 < r_0 \}$ and we consider $r_0 = \{0.1, 0.3\}$. Furthermore, $u_i^{(\ell)}$ and $u_i^{(d)}$ represent respectively the life and death states of the reactants and they are defined in \cite[Chapter 2.4]{daub2013}. Note that we have introduced the $\sin$ function in equations one and three to also test the strength of the procedure on non-polynomial nonlinearities.

For the experimental setup we discretize \ref{cell3d} with $n = 150$ nodes in each of the spatial directions. Therefore, if ${\bf T} \in \RR^{n \times n}$ is defined as in Example~\ref{ex:system}, this yields a system of the form \ref{arraybigsystem}, with 
$$
{\bf A}_{1i} = -c_4{\bf I}_n - \frac{\delta_i}{\ell_x^2}{\bf T}, \quad {\bf A}_{2i} = \frac{\delta_i}{\ell_y^2}{\bf T}, \quad {\bf A}_{3i} = \frac{\delta_i}{\ell_z^2}{\bf T}, \qquad i = 1,2,3,4,
$$
and $\ell_x = \ell_y = \ell_z = 1/(n-1)$, where the additional linear terms have been incorporated into the matrices ${\bf A}_{1i}$. Furthermore, the nonlinear functions ${\mathcal{F}}_1(\pmb{\mathcal{ U}}_1,\pmb{\mathcal{ U}}_2, \pmb{\mathcal{ U}}_3, \pmb{\mathcal{ U}}_4, t)$ and ${\mathcal{F}}_2(\pmb{\mathcal{ U}}_1,\pmb{\mathcal{ U}}_2, \pmb{\mathcal{ U}}_3, \pmb{\mathcal{ U}}_4, t)$ stem from respectively evaluating the nonlinear terms $c_1\sin(u_3u_2)$ and $c_2u_4u_1^3$ elementwise. Notice, furthermore, that ${\mathcal{F}}_3 = -{\mathcal{F}}_1$ and ${\mathcal{F}}_4 = -{\mathcal{F}}_2$, so that only two \MODEIM\ bases are required, instead of four. The constant matrices stemming from the discretization of the constant $c_3$ are treated separately. Finally we consider $\tau = 10^{-2}$ and $n_s = 30$ equispaced snapshots of each $\pmb{\mathcal{ U}}_i$, ${\mathcal{F}}_1$ and ${\mathcal{F}}_2$ in the timespan.

In Table~\ref{celltable} we report, for all four equations and both values of $r_0$, the dimension of the \MPOD\ and \MODEIM\ bases, the online memory requirements ($\#$ of vectors times the length, as before), the online time for $n_{\mathfrak t} = 300$ (for each equation separately) and the average relative error \ref{errormeasure2}.
\begin{table}[htb!]
\caption{\color{black}Example~\ref{ex:cells3d}. Dim. of \MPOD\ and \MODEIM\ bases and further computational detalis for $\tau = 10^{-2}$ and $n = 150$.  \label{celltable}}
\centering
\begin{tabular}{|l|l|r|r|r|r|r|}
\hline             
   &      & \sc pod dim.          &  \sc deim dim.  & \sc online      & \sc online &  \\
   $r_0$ & ${\bf U}_i$     &    {\footnotesize ($k_1/k_2/k_3$)}    &   {\footnotesize ($p_1/p_2/p_3$)}       & \sc memory      & \sc time (s)  & \sc error     \\ \hline
   \multirow{4}{*}{{$0.1$}} & $\BU_1$  & 2/2/2         & 5/5/5  & $21n$      & 1.29 & $3\cdot10^{-4} $      \\ 
 & $\BU_2$  & 2/2/2         & 3/3/3 & $15n$        & 1.20    & $4\cdot10^{-4}$      \\ 
 & $\BU_3$  & 18/18/18         & --  & $54n$       & 1.50    & $3\cdot10^{-4}$       \\ 
 & $\BU_4$  & 9/9/9          & --  & $27n$       & 0.63    & $1\cdot10^{-2} $     \\ \hline
\multirow{4}{*}{{$0.3$}}  & $\BU_1$  & 8/8/8         & 9/9/9  & $51n$       & 1.56    & $3\cdot10^{-4}$       \\ 
 & $\BU_2$  &8/8/8         & 9/9/9  & $51n$      & 1.13    & $3\cdot10^{-4}$     \\ 
 & $\BU_3$  & 43/43/42        & --  & $128n$      & 11.25   & $3\cdot10^{-3}$       \\ 
 & $\BU_4$  & 31/31/30         & --  & $92n$      & 4.27     & $5\cdot10^{-3} $    \\ \hline
 \end{tabular}
\end{table}

We observe a large decrease in the state dimension for both values of $r_0$, with a very acceptable average relative error in all equations. Equations three and four require a larger basis than one and two, but in turn they do not require the additional cost online of \MODEIM\ interpolation and the evaluation of the nonlinear function. Furthermore, we observe that all equations can be solved in a rapid online phase, whereas the full order model needs approximately 4064 seconds to be integrated at $n_{\mathfrak t} = 300$ timesteps, independent of $r_0$. $\square$

}

\end{example}}

\section{Conclusion}
\label{conclusion}
In this paper we have illustrated that systems of the form \ref{systempde}, with linear operators with separable coefficients, discretized by a tensor basis on certain domains, can be treated directly in array form. In this setting, we have extended the \podeim\ model order reduction method to the multilinear setting and illustrated how it can be used to massively reduce the dimension and complexity of systems of ODEs in two and three spatial dimensions. Some very encouraging numerical experiments on difficult problems such as the 2D and 3D viscous Burgers equation, indicate a dramatic decrease in both CPU time and memory requirements in the offline phase for constructing the bases, especially when $d = 3$.

 Nevertheless, the dense Kronecker structure of the reduced order model obtained by the \MPDEIM\ projection would incur unnecessary computational costs in the online phase when $d = 3$. To this end we have shown how the novel {\sc t3-sylv} linear system solver from \cite{simoncini2020} can exploit the structure of the reduced order model, resulting in a large decrease in computational time in the online phase as well. 

Future work would entail an analysis of the number of snapshots required by {\sc ho-pod-deim} in comparison to {\sc pod-deim}. It could also be of interest to extend the {\sc t3-sylv} solver to higher dimensions so that the {\sc ho-pod-deim} strategy can be applied to PDEs with $d > 3$. {\color{black} Furthermore, the extension to the parameter dependent setting can also be considered. This will result in an additional dimension that needs to be treated.} {\color{black} Finally, the presented algorithm can certainly also benefit from a dynamic implementation as presented in the companion manuscript \cite{kirsten2020}.  } 

\section*{Acknowledgments}
We thank Valeria Simoncini for her support and careful reading of earlier versions of this manuscript. {\color{black} We are also grateful to the two anonymous referees for their careful reading and helpful suggestions, which helped improve the presentation.} 

\providecommand{\href}[2]{#2}
\providecommand{\arxiv}[1]{\href{http://arxiv.org/abs/#1}{arXiv:#1}}
\providecommand{\url}[1]{\texttt{#1}}
\providecommand{\urlprefix}{URL }

\medskip
Received xxxx 20xx; revised xxxx 20xx.
\medskip


\begin{thebibliography}{10}

\bibitem{ABG.20}
\newblock A.~Antoulas, C.~Beattie and S.~Gugercin,
\newblock \emph{Interpolatory methods for model reduction},
\newblock SIAM, Philidelphia, 2020.

\bibitem{Ascher1995}
\newblock U.~M. Ascher, S.~J. Ruuth and B.~T. Wetton,
\newblock Implicit-explicit methods for time-dependent partial differential
  equations,
\newblock \emph{SIAM J. Numer. Anal.}, \textbf{32} (1995), 797--823.

\bibitem{astrid2008}
\newblock P.~Astrid, S.~Weiland, K.~Willcox and T.~Backx,
\newblock Missing point estimation in models described by proper orthogonal
  decomposition,
\newblock \emph{IEEE Trans. Autom. Control}, \textbf{53} (2008), 2237--2251.

\bibitem{barrault2004}
\newblock M.~Barrault, Y.~Maday, N.~C. Nguyen and A.~T. Patera,
\newblock An ‘empirical interpolation’ method: application to efficient
  reduced-basis discretization of partial differential equations,
\newblock \emph{C. R. Math. Acad. Sci. Paris}, \textbf{339} (2004), 667--672.

\bibitem{Benner.05}
\newblock P.~Benner, V.~Mehrmann and D.~Sorensen,
\newblock \emph{Dimension reduction of large-scale systems},
\newblock Springer-Verlag, Berlin/Heidelberg, Germany, 2005.

\bibitem{benner2015review}
\newblock P.~Benner, S.~Gugercin and K.~Willcox,
\newblock A survey of projection-based model reduction methods for parametric
  dynamical systems,
\newblock \emph{SIAM Rev}, \textbf{57} (2015), 483--531.

\bibitem{bonomi2017}
\newblock D.~Bonomi, A.~Manzoni and A.~Quarteroni,
\newblock A matrix {DEIM} technique for model reduction of nonlinear
  parametrized problems in cardiac mechanics,
\newblock \emph{Comput. Methods Appl. Mech. Eng.}, \textbf{324} (2017),
  300--326.

\bibitem{chaturantabut2010nonlinear}
\newblock S.~Chaturantabut and D.~C. Sorensen,
\newblock Nonlinear model reduction via discrete empirical interpolation,
\newblock \emph{SIAM J. Sci. Comput.}, \textbf{32} (2010), 2737--2764.

\bibitem{chat2011}
\newblock S.~Chaturantabut and D.~C. Sorensen,
\newblock Application of {POD} and {DEIM} on dimension reduction of non-linear
  miscible viscous fingering in porous media,
\newblock \emph{Math. Comput. Modell. Dyn. Syst.}, \textbf{17} (2011),
  337--353.

\bibitem{daub2012}
\newblock M.~Daub, S.~Waldherr, F.~Allg{\"o}wer, P.~Scheurich and G.~Schneider,
\newblock Death wins against life in a spatially extended apoptosis model,
\newblock \emph{Biosystems}, \textbf{108} (2012), 45--51.

\bibitem{daub2013}
\newblock M.~Daub,
\newblock \emph{Mathematical modeling and numerical simulations of the
  extrinsic pro-apoptotic signaling pathway},
\newblock PhD thesis, University of Stuttgart, 2013.

\bibitem{Autilia2019matri}
\newblock M.~C. D'Autilia, I.~Sgura and V.~Simoncini,
\newblock Matrix-oriented discretization methods for reaction--diffusion
  {PDE}s: {C}omparisons and applications,
\newblock \emph{Computers \& Mathematics with Applications}, 2067--2085.

\bibitem{dewit99}
\newblock A.~De~Wit,
\newblock \emph{Spatial Patterns and Spatiotemporal Dynamics in Chemical
  Systems}, 435--513,
\newblock John Wiley \& Sons, Ltd, 1999,
\newblock
  \urlprefix\url{https://onlinelibrary.wiley.com/doi/abs/10.1002/9780470141687.ch5}.

\bibitem{gugercin2018}
\newblock Z.~Drma\v{c} and S.~Gugercin,
\newblock A new selection operator for the discrete empirical interpolation
  method---improved a priori error bound and extensions,
\newblock \emph{SIAM J. Sci. Comput.}, \textbf{38} (2016), A631--A648.

\bibitem{fletcher1983}
\newblock C.~A. Fletcher,
\newblock Generating exact solutions of the two-dimensional {B}urgers'
  equations,
\newblock \emph{Int. J. Numer. Methods Fluids}, \textbf{3} (1983), 213--216.

\bibitem{gambino2019}
\newblock G.~Gambino, M.~Lombardo and M.~Sammartino,
\newblock Pattern selection in the 2{D} {F}itz{H}ugh--{N}agumo model,
\newblock \emph{Ricerche di Matematica}, \textbf{68} (2019), 535--549.

\bibitem{GAO2017}
\newblock Q.~Gao and M.~Zou,
\newblock An analytical solution for two and three dimensional nonlinear
  {B}urgers' equation,
\newblock \emph{Appl. Math. Modell.}, \textbf{45} (2017), 255 -- 270,
\newblock
  \urlprefix\url{http://www.sciencedirect.com/science/article/pii/S0307904X16306710}.

\bibitem{george2013}
\newblock U.~Z. George, A.~St{\'e}phanou and A.~Madzvamuse,
\newblock Mathematical modelling and numerical simulations of actin dynamics in
  the eukaryotic cell,
\newblock \emph{J Math Biol}, \textbf{66} (2013), 547--593.

\bibitem{golub13}
\newblock G.~H. Golub and C.~F. van Loan,
\newblock \emph{Matrix Computations},
\newblock 4th edition,
\newblock Johns Hopkins University Press, Baltimore, 2013,
\newblock
  \urlprefix\url{http://www.cs.cornell.edu/cv/GVL4/golubandvanloan.htm}.

\bibitem{gu2011}
\newblock C.~Gu,
\newblock {QLMOR}: {A} projection-based nonlinear model order reduction
  approach using quadratic-linear representation of nonlinear systems,
\newblock \emph{IEEE Trans. Comput.-Aided Design Integr. Circuits Syst.},
  \textbf{30} (2011), 1307--1320.

\bibitem{halko2011}
\newblock N.~Halko, P.-G. Martinsson and J.~A. Tropp,
\newblock Finding structure with randomness: {P}robabilistic algorithms for
  constructing approximate matrix decompositions,
\newblock \emph{SIAM Rev}, \textbf{53} (2011), 217--288.

\bibitem{hinze2005}
\newblock M.~Hinze and S.~Volkwein,
\newblock Proper orthogonal decomposition surrogate models for nonlinear
  dynamical systems: {E}rror estimates and suboptimal control,
\newblock in \emph{Dimension reduction of large-scale systems},
\newblock Springer, 2005,
\newblock 261--306.

\bibitem{hodgkin1952}
\newblock A.~L. Hodgkin and A.~F. Huxley,
\newblock A quantitative description of membrane current and its application to
  conduction and excitation in nerve,
\newblock \emph{The Journal of physiology}, \textbf{117} (1952), 500--544.

\bibitem{hundsdorfer2013}
\newblock W.~Hundsdorfer and J.~G. Verwer,
\newblock \emph{Numerical solution of time-dependent
  advection-diffusion-reaction equations}, vol.~33,
\newblock Springer Science \& Business Media, 2013.

\bibitem{karasozen2016}
\newblock B.~Karas{\"o}zen, M.~Uzunca and T.~K{\"u}{\c{c}}{\"u}kseyhan,
\newblock Model order reduction for pattern formation in {F}itzhugh-{N}agumo
  equations,
\newblock in \emph{Numerical Mathematics and Advanced Applications ENUMATH
  2015},
\newblock Springer, 2016,
\newblock 369--377.

\bibitem{karasozen2020}
\newblock B.~Karas{\"o}zen, M.~Uzunca and T.~K{\"u}{\c{c}}{\"u}kseyhan,
\newblock Reduced order optimal control of the convective {F}itzhugh--{N}agumo
  equations,
\newblock \emph{Computers \& Mathematics with Applications}, \textbf{79}
  (2020), 982--995.

\bibitem{karasozen2021}
\newblock B.~Karas{\"o}zen, S.~Y{\i}ld{\i}z and M.~Uzunca,
\newblock Structure preserving model order reduction of shallow water
  equations,
\newblock \emph{Mathematical Methods in the Applied Sciences}, \textbf{44}
  (2021), 476--492.

\bibitem{kirsten2020}
\newblock G.~Kirsten and V.~Simoncini,
\newblock A matrix-oriented {POD-DEIM} algorithm applied to nonlinear
  differential matrix equations,
\newblock \emph{arXiv preprint arXiv:2006.13289}.

\bibitem{kolda2009}
\newblock T.~G. Kolda and B.~W. Bader,
\newblock Tensor decompositions and applications,
\newblock \emph{SIAM Rev}, \textbf{51} (2009), 455--500.

\bibitem{kramer2011}
\newblock B.~Kramer,
\newblock \emph{Model reduction of the coupled {B}urgers equation in
  conservation form},
\newblock PhD thesis, Virginia Tech, 2011.

\bibitem{kramer2019}
\newblock B.~Kramer and K.~E. Willcox,
\newblock Nonlinear model order reduction via lifting transformations and
  proper orthogonal decomposition,
\newblock \emph{AIAA Journal}, \textbf{57} (2019), 2297--2307.

\bibitem{kunisch1999}
\newblock K.~Kunisch and S.~Volkwein,
\newblock Control of the {B}urgers equation by a reduced-order approach using
  proper orthogonal decomposition,
\newblock \emph{J. Optim. theory Appl.}, \textbf{102} (1999), 345--371.

\bibitem{Mainietal.01}
\newblock P.~K. Maini and H.~G. Othmer,
\newblock \emph{Mathematical Models for Biological Pattern Formation},
\newblock The IMA Volumes in Mathematics and its Applications - Frontiers in
  application of Mathematics, Springer-Verlag, New York, 2001.

\bibitem{Malchowetal.08}
\newblock H.~Malchow, S.~Petrovskii and E.~Venturino,
\newblock \emph{Spatiotemporal Patterns in Ecology and Epidemiology: Theory,
  Models, and Simulations},
\newblock Chapman \& Hall, CRC, London, 2008.

\bibitem{matlab2013}
\newblock The MathWorks,
\newblock \emph{MATLAB 7}, r2013b edition, 2013.

\bibitem{minster2020}
\newblock R.~Minster, A.~K. Saibaba and M.~E. Kilmer,
\newblock Randomized algorithms for low-rank tensor decompositions in the
  {T}ucker format,
\newblock \emph{SIAM J. Math. Data Sci.}, \textbf{2} (2020), 189--215.

\bibitem{murray2001}
\newblock J.~Murray,
\newblock \emph{Mathematical biology II: spatial models and biomedical
  applications}, vol.~3,
\newblock Springer-Verlag, 2001.

\bibitem{negri2015}
\newblock F.~Negri, A.~Manzoni and D.~Amsallem,
\newblock Efficient model reduction of parametrized systems by matrix discrete
  empirical interpolation,
\newblock \emph{J. Comput. Phys.}, \textbf{303} (2015), 431--454.

\bibitem{nguyen2008}
\newblock N.-C. Nguyen, A.~T. Patera and J.~Peraire,
\newblock A `best points' interpolation method for efficient approximation of
  parametrized functions,
\newblock \emph{Int J Numer Methods Eng}, \textbf{73} (2008), 521--543.

\bibitem{palitta2016}
\newblock D.~Palitta and V.~Simoncini,
\newblock Matrix-equation-based strategies for convection--diffusion equations,
\newblock \emph{BIT Numerical Mathematics}, \textbf{56} (2016), 751--776.

\bibitem{patera2007reduced}
\newblock A.~T. Patera and G.~Rozza,
\newblock \emph{Reduced basis approximation and a posteriori error estimation
  for parametrized partial differential equations},
\newblock MIT Cambridge, MA, USA, 2007.

\bibitem{Quarteroni.17}
\newblock A.~Quarteroni,
\newblock \emph{Numerical Models for Differential Problems}, vol.~8 of MS\&A -
  Modeling, Simulation and Applications,
\newblock Springer-Verlag, Milan, 2017.

\bibitem{ruuth1995}
\newblock S.~J. Ruuth,
\newblock Implicit-explicit methods for reaction-diffusion problems in pattern
  formation,
\newblock \emph{Journal of Mathematical Biology}, \textbf{34} (1995), 148--176.

\bibitem{sahyoun2014}
\newblock S.~Sahyoun and S.~M. Djouadi,
\newblock Nonlinear model reduction using space vectors clustering {POD} with
  application to the {B}urgers' equation,
\newblock in \emph{2014 American Control Conference},
\newblock IEEE, 2014,
\newblock 1661--1666.

\bibitem{sherratt2001}
\newblock J.~A. Sherratt and M.~A. Chaplain,
\newblock A new mathematical model for avascular tumour growth,
\newblock \emph{J Math Biol}, \textbf{43} (2001), 291--312.

\bibitem{simoncini2020}
\newblock V.~Simoncini,
\newblock Numerical solution of a class of third order tensor linear equations,
\newblock \emph{BUMI}, \textbf{13} (2020), 429--439.

\bibitem{Simoncini2017}
\newblock V.~Simoncini,
\newblock Computational methods for linear matrix equations,
\newblock \emph{SIAM Rev}, \textbf{58} (2016), 377--441.

\bibitem{cstefuanescu2014}
\newblock R.~{\c{S}}tef{\u{a}}nescu, A.~Sandu and I.~M. Navon,
\newblock Comparison of {POD} reduced order strategies for the nonlinear 2{D}
  shallow water equations,
\newblock \emph{Int. J. Numer. Methods Fluids}, \textbf{76} (2014), 497--521.

\bibitem{strikwerda2004}
\newblock J.~C. Strikwerda,
\newblock \emph{Finite difference schemes and partial differential equations},
\newblock SIAM, 2004.

\bibitem{Tveitoetal.10}
\newblock A.~Tveito, H.~P. Langtangen, B.~F. Nielsen and X.~Cai,
\newblock \emph{Elements of Scientific Computing},
\newblock Texts in Computational Science and Engineering, Springer-Verlag,
  Berlin, 2010.

\bibitem{vanag2004}
\newblock V.~K. Vanag,
\newblock Waves and patterns in reaction--diffusion systems.
  {B}elousov--{Z}habotinsky reaction in water-in-oil microemulsions,
\newblock \emph{Phys. Usp.}, \textbf{47} (2004), 923.

\bibitem{vannieuwenhoven2012}
\newblock N.~Vannieuwenhoven, R.~Vandebril and K.~Meerbergen,
\newblock A new truncation strategy for the higher-order singular value
  decomposition,
\newblock \emph{SIAM J. Sci. Comput.}, \textbf{34} (2012), A1027--A1052.

\bibitem{wang2016}
\newblock Y.~Wang, I.~M. Navon, X.~Wang and Y.~Cheng,
\newblock 2{D B}urgers equation with large {R}eynolds number using {POD/DEIM}
  and calibration,
\newblock \emph{Int. J. Numer. Methods Fluids}, \textbf{82} (2016), 909--931.

\end{thebibliography}
\end{document}